\documentclass[10pt,reqno]{amsart} 

\usepackage{amssymb,latexsym}
\usepackage{cite} 

\usepackage[height=190mm,width=130mm]{geometry} 

\usepackage{amsfonts}
\usepackage{amsmath}
\usepackage{amscd}
\usepackage{amssymb}
\usepackage{amsthm}
\usepackage{bbm}
\usepackage{bm}
\usepackage{enumerate}
\usepackage{enumitem}
\usepackage{hyperref}
\usepackage{latexsym}
\usepackage{mathabx}
\usepackage{mathdots}
\usepackage{mathtools}

\usepackage{amsmath}
\usepackage{amssymb}
\usepackage{amsthm}
\usepackage{bbm}
\usepackage{enumerate}
\usepackage{latexsym}
\usepackage{mathdots}
\usepackage{geometry}

\usepackage{hyperref}   

\newcommand{\abar}{\overline{a}}
\newcommand{\bbar}{\overline{b} }

\newcommand{\C}{\mathbb{C}}

\newcommand{\del}{\partial}

\newcommand{\E}{\mathcal{E}}

\newcommand{\F}{\mathcal{F}}
\newcommand{\I}{\mathcal{I}}

\newcommand{\Ip}{\I_{+}}

\newcommand{\SL}{\textup{SL}}

\newcommand{\wt}{\textup{wt}}
\newcommand{\Z}{\mathbb{Z}}

\newcommand{\HH}{\mathbb{H}}

\newcommand{\omt}{\widetilde{\omega}}

\usepackage{amsmath}
\usepackage{amssymb}
\usepackage{amsthm}
\usepackage{bbm}
\usepackage{enumerate}
\usepackage{latexsym}
\usepackage{mathdots}
\usepackage{geometry}

\usepackage{hyperref}   

\theoremstyle{plain}

\newtheorem{proposition}{Proposition}

\theoremstyle{definition}
\newtheorem{definition}{Definition}

\theoremstyle{remark}
\newtheorem{remark}{Remark}

\numberwithin{equation}{section} 

\newcommand{\rmap}{\longrightarrow}

\newcommand{\g}{\ensuremath{\Gamma}}

\newcommand{\ps}{{\raise 1pt\hbox{\tiny (}}}

\newcommand{\pss}{{\raise 1pt\hbox{\tiny [}}}
\newcommand{\pdd}{{\raise 1pt\hbox{\tiny ]}}}
\newcommand{\pd}{{\raise 1pt\hbox{\tiny )}}}

\newcommand{\bs}{{\raise 1pt\hbox{\tiny [}}}
\newcommand{\bd}{{\raise 1pt\hbox{\tiny ]}}}

\def\cross{\mathinner{\mathrel{\raise0.8pt\hbox{$\scriptstyle>$}}
                 \joinrel\mathrel\triangleleft}}

\usepackage{citehack}
\usepackage{amsmath,amsthm}
\usepackage{amsfonts}
\usepackage{amssymb}
\usepackage{longtable}
\usepackage[matrix,arrow,curve]{xy}

\usepackage{lipsum}

\def\V{\mathcal{V}}

\def\K{\mathcal{K}}

\usepackage{stackrel}

\newcommand{\be}{\begin{equation}}
\newcommand{\ee}{\end{equation}}

\newcommand{\nn}{\nonumber \\}

\newcommand{\one}{\mathbf{1}}

\newcommand{\nc}{\newcommand}
\nc{\cali}{\mathcal}
\nc{\on}{\operatorname}
\nc{\Wick}{{\mb :}}

\nc{\ddz}{\frac{\partial}{\partial z}}
\nc{\ch}{\mbox{ch}}
\nc{\Oo}{{\cali O}}
\nc{\cond}{|\,}
\nc{\bib}{\bibitem}
\nc{\pone}{\Pro^1}
\nc{\pa}{\partial}
\nc{\arr}{\rightarrow}
\nc{\larr}{\longrightarrow}
\nc{\ket}{\rangle}
\nc{\bra}{\langle}
\nc{\gam}{\bar{\gamma}}
\nc{\q}{\widetilde{Q}}
\nc{\ep}{\epsilon}
\nc{\su}{\widehat{{\mf s}{\mf l}}_2}
\nc{\sw}{{\mf s}{\mf l}}
\nc{\h}{{\mf h}}
\nc{\n}{{\mf n}}
\nc{\ab}{\mf{a}}
\nc{\is}{{\mb i}}
\nc{\js}{{\mb j}}
\nc{\bi}{\bibitem}
\nc{\He}{{\cali H}}
\nc{\inv}{^{-1}}
\nc{\ol}{\overline}
\nc{\wh}{\widehat}
\nc{\dst}{\displaystyle}

\nc{\delt}{\partial_t}
\nc{\ddt}{\frac{\partial}{\partial t}}
\nc{\delx}{\partial_x}
\nc{\mb}{\mathbf}
\nc{\mf}{\mathfrak}

\nc{\mbb}{\mathbb}
\nc{\Ctt}{\C((t))}
\nc{\Ct}{\C[t,t\inv]}

\nc{\un}{\underline}
\nc{\mc}{\mathcal}
\nc{\BB}{{\mc B}}
\nc{\bb}{{\mf b}}
\nc{\kk}{{\mf k}}
\nc{\frob}{\times}
\nc{\sm}{\setminus}
\nc{\Pp}{{\mathbb P}^1}
\nc{\Aa}{{\mc A}}

\nc{\AutO}{\on{Aut}\Oo}
\nc{\AUTO}{\un{\on{Aut}}\Oo}
\nc{\AUTK}{\un{\on{Aut}}\K}
\nc{\Heout}{\He_{\out}}
\nc{\Hetil}{{\widetilde\He}}
\nc{\wb}{\overline}

\nc{\Res}{\on{Res}}
\nc{\pitil}{\Pi}
\nc{\Ctil}{\wt{C}}
\nc{\auto}{\on{Aut} \Oo}
\nc{\phitil}{\wt{\phi}}
\nc{\gz}{\g_{\vec z}}
\nc{\tensorM}{\bigotimes_{i=1}^N{\mathbb M}_i}
\nc{\tensorW}{\bigotimes_{i=1}^N W_{\nu_i,k}}
\nc{\out}{\on{out}}

\nc{\m}{{\mathfrak m}}

\nc{\gx}{\g^0_{\vec x}}

\nc{\hx}{\He^0_{\vec x}}
\nc{\tensorpi}{\pi_{\nu_1,\ldots,\nu_N}^\kappa}
\nc{\Phizw}{\Phi_{\vec w}({\vec z})}
\nc{\Pro}{{\mathbb P}}

\nc{\De}{\Delta}

\nc{\us}{\underset}

\nc{\Ll}{\mc L}
\nc{\dR}{\on{dR}}

\nc{\T}{{\mc T}}

\nc{\Xn}{\overset{\circ}X{}^n} \nc{\Dn}{\overset{\circ}D{}^n}
\nc{\Dxn}{\overset{\circ}D{}^n_x} \nc{\varphitil}{\wt{\varphi}}

\nc{\lf}{{\mf l}}
\nc{\GL}{{}^L G}
\nc{\Vir}{\on{Vir}}

\begin{document}
\title[Reduction cohomology of Riemann surfaces]
{Reduction cohomology of Riemann surfaces}  
\author{A. Zuevsky}
\address{Institute of Mathematics \\ Czech Academy of Sciences\\ Zitna 25, Prague \\ Czech Republic}

\email{zuevsky@yahoo.com}







\begin{abstract}
 We study the algebraic conditions leading to the chain property of complexes for 
vertex operator algebra $n$-point functions (with their convergence assumed)
 with differential being defined through reduction formulas. 
The notion of the reduction cohomology of Riemann surfaces is introduced.  
Algebraic, geometrical, and cohomological meanings of reduction formulas is clarified. 
A counterpart of the Bott-Segal theorem for Riemann surfaces in terms of the reductions cohomology is proven. 
It is shown that the reduction cohomology is given by the cohomology of $n$-point 
connections over the vertex operator algebra bundle defined on a genus $g$ Riemann surface $\Sigma^{(g)}$. 
The reduction cohomology for a vertex operator algebra with formal parameters identified 
with local coordinates around marked points on $\Sigma^{(g)}$ is found  
in terms of the space of analytical continuations of solutions to Knizhnik-Zamolodchikov equations. 
For the reduction cohomology, 
the Euler-Poincare formula is derived. 
Examples for various genera and vertex operator cluster algebras are provided.  
\end{abstract}

\keywords{Cohomology; Vertex algebras; Riemann surfaces; Cluster algebras}

\maketitle

\section{Introduction} 
The natural problem of computing continuous cohomologies 
for non-commutative structures  
on manifolds
 has proven to be a subject of great geometrical interest \cite{BS, Kaw, PT, Fei, Fuks, Wag}.   
As it was demonstrated in \cite{Fei, Wag}, the ordinary Gelfand-Fuks cohomology of the Lie algebra of 
holomorphic vector fields on complex manifolds 
 turns to be not the most effective and general one.  
For Riemann surfaces, and even for higher dimension complex manifolds, the classical cohomology of 
vector fields becomes trivial \cite{Kaw}.   
The Lie algebra of 
holomorphic vector fields is not always an interesting Lie algebra. For example, it is zero for a 
compact Riemann surface of genus greater than one, and one looks for other algebraic objects having locally 
the same cohomology.   
In \cite{Fei} Feigin obtained various results concerning (co)-homology of  
 the Lie algebra cosimplicial objects of holomorphic vector fields $Lie(M)$.   
Despite results in previous approaches, it is desirable to     
  find a way to enrich cohomological structures which motivates 
 constructions of more refined cohomology description for non-commutative 
algebraic structures.   
In the seminal paper \cite{BS}, the authors  
proved that the Gelfand-Fuks cohomology $H^*(Vect(M))$ of vector fields on a smooth compact manifold $M$ 
 is isomorphic to the singular cohomology of the space of continuous cross sections of a certain 
fibre bundle over $M$. 
An important problem of revealing relations between non-commutative structures 
and geometrical objects on complex manifolds 
still remains underinvestigeted in the literature \cite{PT}.  

The main idea of this paper is to introduce and compute the reduction cohomology vertex operator algebras 
\cite{B, FLM, FHL, K}   
with formal parameters considered 
as local coordinates on a genus $g$ compact Riemann surface \cite{FK, Bo, Gu, A} is to study cohomology of  
spaces of converging functions with respect to 
adding new sets of pairs of vertex operator algebra elements and corresponding formal parameters 
\cite{BZF, DGM, DL, EO, FMS}.   
Due to structure of correlations functions \cite{FMS} and reduction relations 
\cite{Y, Zhu, MTZ, GT, TW} among them, 
one can form chain complexes of $n$-point functions (with their convergence assumed) that 
 are fine enough to describe local geometry 
of Riemann surfaces.  
Another meaning of the reduction cohomology of Riemann surfaces is how 
sections of certain bundles with values in higher genus generalizations of elliptic funtions change 
along growing of number of marked points on a surface. 
In addition to that we prove a version of the Bott-Segal \cite{BS, Wag, PT}  
theorem for compact Riemann surfaces of arbitrary genus. 
It relates the reduction cohomology with cohomology of a space of sections of certain vertex operator algebra bundle 
\cite{BZF}. 
Our new vertex operator algebra approach to cohomology of compact Riemann surfaces    
involves Lie-algebraic formal series with applications of techniques used in surgery of spheres \cite{Huang}.   
In contrast to more geometrical methods, e.g.,  in ordinary cosimplicial cohomology for 
Lie algebras \cite{Fei, Wag},     
the reduction cohomology pays more attention to the analytical structure of elements 
 of chain complex spaces. 
Computational methods involving reduction formulas proved  their effectiveness in  
 conformal field theory, geometrical descriptions of intertwined modules for 
Lie algebras, 
 and differential geometry of integrable models.    
In section \ref{chain} we give the definition of the reduction cohomology and prove a proposition 
relating it to cohomology of a vertex operator algebra bundle in terms of $n$-point connections.
In section \ref{cohomology} the main proposition expressing the reduction cohomology in terms of 
spaces of auxiliary functions on Riemann surfaces is proven. 
Results  of this paper are useful for cosimplisial cohomology theory of smooth manifolds, 
  generalizations of the Bott-Segal theorem
have their consequences in conformal field theory 
\cite{Fei, Wag}, deformation theory \cite{O}, non-commutative geometry, modular forms, 
 and the theory of foliations.   
\section{Chain complex for vertex operator algebra $n$-point functions} 
\label{chain}
\subsection{Spaces of $n$-point correlation functions}
In this section we recall definitions and some properties of 
 correlation functions for vertex operator algebras on Riemann surfaces \cite{FHL, Zhu, FMS, DGM}.  
Let us fix a vertex operator algebra $V$. Depending on its commutation relations and configuration of a genus $g$ 
Riemann surface $\Sigma^{(g)}$, the space of all $V$ $n$-point functions 
 can contain various forms of complex functions defined on $\Sigma^{(g)}$.   

We denote by ${\bf v}_n=(v_1, \ldots, v_n) \in V^{\otimes n}$ a tuple of vertex operator algebra elements (see 
Appendix \ref{vertex} for definition of a vertex operator algebra). 
Mark $n$ points ${\bf p}_n=(p_1, \ldots, p_n)$ on a Riemann surface of genus $\Sigma^{(g)}$. 
Denote by ${\bf z}_n=(z_1, \ldots, z_n)$ local coordinates around ${\bf p}_n \in \Sigma^{(g)}$. 
Let us introduce the notation: ${\bf x}_n= \left({\bf v}_n, {\bf z}_n \right)$.
On a genus $g$ Riemann surfaces an $n \ge 0$-point correlation function 
$\F^{(g)}_n \left({\bf x}_n, B^{(g)}\right)$ (here $B^{(g)}$ denotes the set of its parameters) has 
certain specific form depending on $g$, construction of a Riemann surface $\Sigma^{(g)}$, type of 
conformal field theory 
 model used for definitions of $n$-point functions, 
 and the type of commutation relations for $V$-elements. 
In addition to that, it depends on 
 a set of moduli parameters $B^{(g)} \in {\mathcal B}^{(g)}$  
 where we denote by ${\mathcal B}^{(g)}$ a subset of   
 the moduli space of genus $g$ Riemann surfaces $\Sigma^{(g)}$ obtained by specific ways of sewing of lower 
genus Riemann surfaces.  
In particular, for $n$-point functions considered on Riemann surfaces, $B^{(g)}$ is an element of the space  \cite{Y}  
 ${\mathcal B}^{(g)}= \left(\Omega^{(i)}; (\epsilon_i);  (\rho_i); (\mu_i) \right)$,  
$0 \le i \le g$, where vectors $(\epsilon_i)$, $(\rho_i)$, are Riemann surface sewing parameters \cite{Y}, 
 $(\mu_i)$ are any further modular space parameters,     
and $\Omega^{(i)}$ 
is the period matrix of corresponding Riemann surface $\Sigma^{(g_i)}$ used in the procedure of $\Sigma^{(g)}$ 
construction. 
\begin{definition}   
For a Riemann surface $\Sigma$ of genus $g$, a $V$-module $W$, and $n \ge 0$, ${\bf x}_n$ on $\Sigma$,   
  we consider the spaces of $n$-point correlation functions     
\[
C^n_{(g)}(W)= \left\{ 
\F_{W, n}^{(g)} \left({\bf x}_n, B^{(g)}\right)
\right\}.   
\]
\end{definition}
Note that we choose the same $V$-module $W$ for all ${\bf x}_n$.  
In what follows, we will omit where possible $W$ and $B^{(g)}$.
 The co-boundary operator $\delta^n_{(g)}(v_{n+1})$ on $C_{(g)}^n(W)$-space is defined     
 according to the reduction formulas for $V$-module $W$ on a genus $g$ Riemann surface 
(cf. particular examples in subsections \ref{sphere}--\ref{cluster}, \cite{Zhu, MTZ, GT, TW}). 
\begin{definition}
 For $g \ge 0$, $n \ge 0$, and any $x_{n+1} \in V\times \C$, define     
\begin{eqnarray}
\label{delta_operator}
\delta_{(g)}^n: C_{(g)}^n(W) & {\rightarrow }& C_{(g)}^{n+1}(W),     
\nn
\delta^{(g)}_n =  H^{(g)} &=& H^{(g)}_1  + H^{(g)}_2,    
\end{eqnarray}
with 
operators
$H_1^{(g)} (x_{n+1})$, $H_2^{(g)} (x_{n+1})$ given by 
\begin{eqnarray}
\label{poros}
 && H^{(g)}_1 (x_{n+1}) \; \F_{n}^{(g)} \left( {\bf x}_n \right) = 
\sum\limits_{l=1}^{l(g)} 
  f_1^{(g)} \left( {\bf x}_{n+1}, l \right)  \;  T^{(g)}_{l} \F^{(g)}_{n} \left( {\bf x}_n \right),   
\nn
&& H_2^{(g)} ( x_{n+1} )\; \F_{n}^{(g)} \left( {\bf x}_n \right) 
 =
\nonumber
 \sum\limits_{k=1}^{n} \sum\limits_{m \ge 0}   
  f_2^{(g)} ( {\bf x}_{n+1}, k, m ) 
\nn
&& \qquad \qquad \qquad \qquad \qquad \qquad \qquad \qquad
\cdot  T^{(g)}_k( v_{n+1}(m) )\; \F^{(g)}_{W, n} \left( {\bf x}_n \right),   
\end{eqnarray} 
where $l(g) \ge 0$ is a constant depending on $g$, and the meaning of indexes $1 \le k \le n$, 
$1 \le l \le l(g)$, $m \ge 0$ explained below. 
\end{definition}
 Operator-valued functions $f^{(g)}_1\left( {\bf x}_{n+1}, l \right) \; T^{(g)}_{l}$, 
$f_2^{(g)}({\bf x}_{n+1}, k, m )$. $T^{(g)}_k(v_{n+1}(m))$   
 depend on genus of a Riemann surface $\Sigma^{(g)}$.  
$T^{(g)}_{l}$    
and $T^{(g)}_k(v(m))$ are operators 
of insertion of certain function of
vertex operator algebra modes into $\F^{(g)}_{W, n}\left( {\bf x}_n \right)$ at 
the $k$-th entry: 
\begin{eqnarray}
T^{(g)}_{l, k}(v(m))\; \F^{(g)}_{W, n}\left( {\bf x}_n \right) &=&  
\F^{(g)}_{W, n} \left( \ldots,  \left(T^{(g)}_{l} (v(m)). x_{n}\right)_k, \ldots \right), 
\nonumber
\end{eqnarray}
where we use the notation 
\[
(\gamma.)_k\; x_{n} =  
\left(x_1,  \ldots,  \gamma.x_k,  \ldots, x_n \right),  
\]
for an operator $\gamma$ acting on $k$-th entry.  
 Note that commutation properties of $H_1^{(g)}(x_{n+1})$ and $H_2^{(g)}(x_{n+1})$ depend on genus $g$. 
The reduction formulas have the form:  
\begin{equation}
\label{reduction}
\F^{(g)}_{W, n+1}\left( {\bf x}_{n+1} \right)= H^{(g)}(x_{n+1}) \; 
\F^{(g)}_{W, n}\left( {\bf x}_n \right), 
\end{equation}
for $n \ge 0$.  
\begin{remark}
The author's conjecture (based on \cite{Y, MT, MTZ, TZ, GT, TW, BKT}) is that for all possible
 configurations of a genus $g$  
Riemann surface, the form of the reduction relations \eqref{poros} is coinvariant and it is given by a sum of 
operators acting on the an $n$-point correlation function  
\[
\F^{(g)}_{W, n} \left({\bf x}_{n+1} \right)=\sum\limits_{1 \le l \le l(g) \atop m \ge 0, k \ge 0}
 f^{(g)}\left({\bf x}_{n+1}, k,l,m \right) \; T^{(g)}( k, l, m)\; 
\F^{(g)}_{W, n} \left({\bf x}_{n}  \right).  
\]
\end{remark}
\begin{remark}
The reductions formulas have an interpretation in terms of torsors \cite{BZF} (Chapter 6). In such formulation
${\bf x}_n$ is a torsor with respect to the group of transformation of the space $V^{\otimes \; n}\times \C^n$. 
In particular, from \eqref{reduction} we see that $T^{(g)}_k(u(m))$-operators act on $V^{\otimes \; n}$-entries 
of ${\bf x}_n$ in $\F^{(g)}_{W, n}({\bf x}_n, B^{(g)})$, while $T^{(g)}_{k,l}$-operators act on 
${\bf z}_n$ of $\F^{(g)}_{W, n}({\bf x}_n, B^{(g)})$ 
as a complex function. 
\end{remark}
For $n \ge 0$, 
let us denote by ${\mathfrak V}_n$ the subsets of all ${\bf x}_{n+1} \in V^{\otimes (n+1)}\times \C^{n+1}$,
 such that  
 the chain condition 
\begin{equation}
\label{torba}
H^{(g)}(x_{n+1})\; H^{(g)}(x_{n}) \; \F^{(g)}_{W, n}\left({\bf x}_n \right)=0,   
\end{equation}
 for the coboundary operators \eqref{poros} for complexes $C^n_{(g)}(W)$   
 is satisfied.   

Explicitly, 
the chain condition \eqref{torba} 
leads to an infinite $n \ge 0$ set of equations involving functions 
$f_1^{(g)}\left({\bf x}_{n+1}, l \right)$,   
$f_2^{(g)}\left({\bf x}_{n+1}, k, m \right)$, and   
 $\F^{(g)}_{W, n}\left({\bf x}_n \right)$:  
\begin{eqnarray}
\label{conditions}
&&\Big(\sum\limits_{l'=1, \; l=1}^{l(g)}   
  f_1^{(g)}\left({\bf x}_{n+2}, l' \right)  
\; f_1^{(g)}\left({\bf x}_{n+1}, l;  B^{(g)}\right) \; T^{(g)}_{l'} \; T^{(g)}_l    
\nn
&&+ \sum\limits_{l'=1, }^{l(g)}  \sum\limits_{k=1}^{n} \sum\limits_{m \ge 0} 
 f^{(g)}_1\left({\bf x}_{n+2}, l'  \right)   
  f_2^{(g)}( m, k ) \;  
\;  T^{(g)}_{l'} \; T^{(g)}_k(v_{n+1}(m)) 
\nn
&&+ \sum\limits_{k'=1}^{n+1} \sum\limits_{m' \ge 0} 
 \sum\limits_{l=1, }^{l(g)}  
  f_2^{(g)}\left( k', m' \right) f_1^{(g)}\left({\bf x}_{n+2}, l  \right)   
   \;T^{(g)}_{k'}(v_{n+1}(m')) \; T^{(g)}_l   
\nn
 &&+ \sum\limits_{k'=1}^{n+1} \sum\limits_{k=1}^{n} \; \sum\limits_{m', \atop m \ge 0}   
  f_2^{(g)}\left( k',m' \right) \;   f_2^{(g)}\left( k, m \right)\; 
T^{(g)}_{k'}(v_{n+2}(m'))  \; T^{(g)}_k(v_{n+1}(m)) \Big)\; 
\nn
&& \qquad  .\F^{(g)}_{W, n}\left( {\bf x}_n \right)=0.
\end{eqnarray}
\begin{remark}
\eqref{conditions} contain finite series and    
 narrows the space of compatible $n$-point functions.   
The subspaces of $C_{(g)}^n(W)$, $g \ge 0$, $n \ge 0$, of genus $g$ $n$-point functions 
such that the condition \eqref{conditions} is fulfiled 
for reduction cohomology complexes 
are non-empty. 
%
Indeed, for all $g$,  
 the condition \eqref{conditions} represents an infinite $n \ge 0$ set of functional-differential 
equations (with finite number of summands) on converging complex functions 
$\F^{(g)}_{W, n}({\bf x}_n )$ defined for  
 $n$ local complex variables on a Riemann surface 
of genus $g$ with functional coefficients $f^{(g)}_1\left(v_{n+1}, l \right)$ 
and $f^{(g)}_2 \left(k, m \right)$ 
 (in our examples in subsection \ref{sphere}--\ref{cluster}, these are generalizations of genus $g$ elliptic 
 functions)
 on $\Sigma^{(g)}$. 
Note that (see examples in Sectiion \ref{examples}), all vertex operator algebra elements of ${\bf v}_n\in V^{\otimes n}$,  
 as non-commutative parameters are not present in final form of functional-differential equations 
since they incorporated into either matrix elements, traces, etc. 
According to the theory of such equations \cite{FK, Gu},  
each equation in the infinite set of \eqref{conditions} always have a solution in domains 
they are defined. 
Thus, there always exist solutions of \eqref{conditions} defining $\F^{(g)}_{W, n} \in C^n_{(g)}(W)$,
 and they are not empty.  
\end{remark}
\begin{definition}
The spaces with conditions \eqref{conditions} constitute a semi-infinite chain complex
\begin{equation}
\label{buzova}
0 \rmap 
C^0_{(g)} \stackrel {\delta^{0}_{(g)}} {\rmap} C^1_{(g)}  
\stackrel {\delta^{1}_{(g)}} {\rmap} 
\ldots \rmap C^{n-1}_{(g)} \stackrel{\delta^{n-1}}{\rmap} C^n_{(g)}\rmap \ldots. 
\end{equation} 
For $n \ge 1$, 
we call corresponding cohomology 
\begin{equation}
\label{pusto}
H^n_{(g)}={\rm Ker}\; \delta^{n}_{(g)}/{\rm Im}\; \delta^{n-1}_{(g)},    
\end{equation}
the $n$-th reduction cohomology of a vertex operator algebra $V$-module $W$ 
on a compact Riemann surface $\Sigma^{(g)}$ of genus $g$.  
\end{definition}
In particular, 
the operators $T^{(g)}_{l}$, $0 \le l \le l(g)$, $T^{(g)}_k(u(m))$, $m\ge 0$, $1 \le k \le n$, 
 form a set of generators of an infinite-dimensional continual Lie algebra $\mathfrak g(V)$ 
 endowed with a natural grading indexed by $l$ and $m$.    

Indeed, we set the space of functions $\F^{(g)}_{W, n}$
   as the base algebra \cite{Sav} for the continual Lie algebra $\mathfrak g(V)$, and 
 the generators as 
\begin{eqnarray}
X_{0, l} \left(\F_{W, n}^{(g)}\left({\bf x}_n \right)\right) &=&   
T^{(g)}_{l}\left(\F_{W, n}^{(g)} \left({\bf x}_n \right)\right), 
\nn
X_{k, m} \left(\F_{W, n}^{(g)}\left({\bf x}_n \right)\right) &=&    
T^{(g)}_k(u(m))\left( \F_{W, n}^{(g)} \left({\bf x}_n  \right)\right). 
\end{eqnarray}
for $0 \le l \le l(g)$, $m \ge 0$, $1 \le k \le n$.    
Then the commutation relations for vertex operator algebra modes $v_{n+1}(m).$ in the action of 
operators $T^{(g)}_l$ and $T^{(g)}_{k, l}$
 on ${v}_k$, $1 \le k \le n$ inside $\F^{(g)}_{W, n}$ represent the 
commutation relations of the continual Lie algebra $\mathfrak g(V)$.  
Jacobi identities for $\mathfrak g(V)$ follow from Jacobi identities 
\eqref{vertex operator algebraJac} for a vertex operator algebra $V$. 
\begin{remark}
Recall that we consider genus $g$ Riemann surfaces resulting from combinations of 
sewing procedures of \cite{Y}. Accordingly, due to \cite{MT, TZ, TZ1, TZ2, GT, TW}, corresponding genus $g$  
 $n$-point functions are obtained coherently by combining lower genus functions.  
Then, relations among $n$-point functions of various genera appear. 
One is able to consider a cohomology theory taking into account such relations. 
For instance,  for the $\epsilon$-formalism of \cite{Y} one has 
\begin{eqnarray}
C^{(g)}_n(W) &\rightarrow& C^{(g_1)}_{n_1+1}(W_1) \times C^{(g_2)}_{n_2+1}(W_2), 
\\
\F^{(g)}_{W, n}\left({\bf x}_n \right)  
\nonumber
&=& \sum\limits_{l \ge 0} \epsilon^l \F^{(g_1)}_{W, n+1} \left( (\bar{u}, z), {\bf x}_n; B^{(g_1)} \right) \; 
\F^{(g_2)}_{w, n+1} \left( (u, z), {\bf x}'_n; B^{(g_2)} \right), 
\end{eqnarray}
for $g=g_1+g_2$, $n=n_1 + n_2$, 
and $W=W_1\otimes W_2$.  
In the
$\rho$-formalism \cite{Y} one has  
\begin{eqnarray}
C^{(g)}_n(W) &\rightarrow& C^{(g-1)}_{n_1+1}(W), 
\\
\nonumber
\F^{(g)}_{W, n}\left({\bf x}_n \right) &=& \sum\limits_{l \ge 0} \rho^l \; 
 \F^{(g-1)}_{W, n+2}\left( (u, z), {\bf x}_n; 
 (\bar{u}, z), {\bf x}'_n; B^{(g-1)}\right). 
\end{eqnarray}
\end{remark}
\subsection{Geometrical meaning of reduction formulas and conditions \eqref{conditions}}    
In this section we show that the reduction formulas have the form of multipoint connections   
generalizing ordinary holomorphic connections on complex curves \cite{BZF}.  

\subsubsection{Holomorphic $n$-point connections}
Let us define the notion of a multipoint connection which will be usefull for identifying 
reduction cohomology in section \ref{cohomology}. 
Motivated by the definition of a holomorphic connection  
  for a vertex operator algebra bundle (cf. Section 6, \cite{BZF} and \cite{Gu}) over   
a smooth complex curve, we introduce the definition of 
the multiple point connection over $\Sigma^{(g)}$.  
\begin{definition}
\label{mpconnection}
Let $\V$ be a holomorphic vector bundle over $\Sigma^{(g)}$, and $\mathcal X_0 \subset \Sigma^{(g)}$ be its subdomain.
Denote by ${\mathcal S \mathcal V}$ the space of sections of $\V$.  
A multi-point 
connection $\mathcal G$ on $\V$    
is a $\C$-multi-linear map 
\[
\mathcal G: \left(\Sigma^{(g)}\right){}^{\times n} \times V^{\otimes n}   \to \C,  
\] 
such that for any holomorphic function $f$, and two sections $\phi(p)$ and $\psi(p')$ of $\V$ at points  
$p$ and $p'$ on $\mathcal X_0 \subset \Sigma^{(g)}$ correspondingly, we have 
\begin{equation}
\label{locus}
\sum\limits_{q, q' \in \mathcal X_0 \subset \Sigma^{(g)}} \mathcal G\left( f(\psi(q)).\phi(q') \right) = f(\psi(p')) \; 
\mathcal G \left(\phi(p) \right) + f(\phi(p)) \; \mathcal G\left(\psi(p') \right),    
\end{equation}
where the summation on left hand side is performed over locuses of points $q$, $q'$ on $\mathcal X_0$.   
We denote by ${\mathcal Con}_{n}$ the space of $n$-point connections defined over $\Sigma^{(g)}$. 
\end{definition}
Geometrically, for a vector bundle $\V$ defined over $\Sigma^{(g)}$, 
a multi-point connection \eqref{locus} relates two sections $\phi$ and $\psi$ at points $p$ and $p'$
with a number of sections on $\mathcal X_0 \subset \Sigma^{(g)}$.  
\begin{definition}
We call
\begin{equation}
\label{gform}
G(\phi, \psi) = f(\phi(p)) \; \mathcal G\left(\psi(p') \right)  + f(\psi(p')) \; \mathcal G \left(\phi(p) \right)   
- \sum\limits_{q, q' \in \mathcal X_0 \subset \mathcal X} \mathcal G\left( f(\psi(q')).\phi(q) \right), 
\end{equation}
the form of a $n$-point connection $\mathcal G$. 
The space of $n$-point connection forms will be denoted by $G^n$.   
\end{definition}
Here we prove the following 
\begin{proposition}
\label{pisaka}
$n$-point correlation functions of the space  $\left\{ \F^{(g)}_{W, n} \left({\bf x}_n \right), n \ge 0  \right\}$ 
 form $n$-point connections on the space of sections of the vertex operator algebra bundle $\V$ associated to $V$.   
For $n\ge 0$, the reduction cohomology of a compact Riemann surface of genus $g$ is   
\begin{equation}
\label{chek}
H^n_{(g)}(W)  = H^n_{(g)}(\mathcal S\V)= {\mathcal Con}^{n}/G^{n-1},     
\end{equation}
is isomorphic to the  
cohomology of  
 the space of $\V$-sections.      
\end{proposition}
\begin{remark}
Proposition \ref{pisaka} is a vertex operator algebra version  
of the main proposition of \cite{BS, Wag}, i.e., the Bott-Segal theorem for Riemann surfaces.      
\end{remark}
\begin{proof}
In \cite{BZF} (Chapter 6, subsection 6.5.3) the vertex operator bundle $\V$ was explicitely constructed. 
It is easy to see that $n$-point connections are holomorphic connection on the bundle $\V$ 
with the following identifications. 
For non-vanishing $f(\phi(p))$ let us write \eqref{reduction} as
\begin{eqnarray}
 \mathcal G\left(\psi(p') \right) 
&=& -  \frac{f(\psi(p'))}{ f(\phi(p))} \; 
\mathcal G \left(\phi(p) \right) 
 + \frac{1}{f(\phi(p)) } \sum\limits_{q, q' \in \mathcal X_0 \subset \mathcal X} 
 \mathcal G\left( f(\psi(q)).\phi(q') \right).      
\end{eqnarray}
Let us set  
\begin{eqnarray}
\label{identifications}
\mathcal G &=&\F^{(g)}_{W, n}\left({\bf x}_n \right),
\nn  
\psi(p')&=&\left({\bf x}_{n+1} \right), 
\nn
\phi(p)&=&\left({\bf x_n} \right), 
\nn
 \mathcal G\left( f(\psi(q)).\phi(q') \right)  &=& T^{(g)}_k(v(m))\; \F^{(g)}_{W, n}\left( {\bf x}_n \right), 
\nn
-  \frac{f(\psi(p'))}{ f(\phi(p))} \; 
\mathcal G \left(\phi(p) \right)&=&
\sum\limits_{l=1}^{l(g)} 
  f^{(g)}_1\left(v_{n+1}, l \right)  \;  T^{(g)}_l \F^{(g)}_{W, n}\left( {\bf x}_n \right), 
\nn
 \frac{1}{f(\phi(p)) }  \sum\limits_{{q}_n, { q'}_n \in \atop {\mathcal X}_0  \subset \Sigma^{(g)}} 
\mathcal G\left( f(\psi(q)).\phi(q') \right) 
&=& 
 \sum\limits_{k=1}^{n} \sum\limits_{m \ge 0}  
  f^{(g)}_2\left( k, m \right) 
\nn
&& \qquad \qquad 
\; \cdot T^{(g)}_k(v(m))\; \F^{(g)}_{W, n} \left( {\bf x}_n \right). 
\nn
&&
\end{eqnarray}
Thus, the formula \eqref{identifications} gives \eqref{reduction}. 
Recall \cite{BZF} the construction of the vertex operator algebra bundle $\V$.
According to Proposition 6.5.4 of \cite{BZF}, one canonically (i.e., coordinate independently) associates 
${\rm End} \V$-valued sections $\mathcal Y_p$ of $\V^*$ (the bundle dual to $\V$) 
to matrix elements of a number of vertex operators  
on appropriate punctured disks around points with local coordinates ${\bf z}_n$ on $\Sigma^{(g)}$.  
The spaces of such $\V$-sections for each $n$ are described by identifications \eqref{identifications}. 
Taking into account the construction of Section 6 (subsection 6.6.1, in particular, construction 6.6.4, 
and Proposition 6.6.7) 
 of \cite{BZF}, we see that $n$-point functions  
are connections on the space of sections of $\V$, 
 and the reduction cohomology \eqref{pusto} is represented by \eqref{chek}. 
\end{proof}
For the chain condition \eqref{conditions} we have 
%
\begin{eqnarray*}
 0=\mathcal G\left(\chi(p'') \right) 
&=& -  \frac{f(\chi(p''))}{ f(\psi(p'))} {\mathcal G}(\psi(p'))
 + \frac{1}{f(\psi(p')) } \sum\limits_{{\widetilde {q}}_n, {\widetilde { q}}'_n \in \mathcal X_0 \subset \Sigma^{(g)}}  
 {\mathcal G}\left( f(\chi({\widetilde q})).\psi({\widetilde q}') \right),      
\end{eqnarray*}
\begin{eqnarray}
 0&=&\mathcal G\left(\chi(p'') \right) 
\nn
&=&  \frac{f(\chi(p''))} 
  { f(\phi(p))} \; 
\mathcal G \left(\phi(p) \right) 
  - \frac{f(\chi(p''))} { f(\psi(p')) \; f(\phi(p))}  \sum\limits_{{ q}_n, { q}'_n \in \mathcal X_0 \subset \Sigma^{(g)}} 
 \mathcal G\left( f(\psi(q)).\phi(q') \right)
\nn
 &+& \frac{1}{f(\psi(p')) }  \sum\limits_{{\widetilde {q}}_n, {\widetilde {q}}'_n \in \mathcal X_0 \subset \Sigma^{(g)}} 
{\mathcal G} \left( f(\chi({\widetilde q})).\psi({\widetilde q}') \right),      
\end{eqnarray}
%
The geometrical meaning of \eqref{conditions} consists in the following. 
Since in \eqref{poros} operators act on vertex operator algebra elements only, 
 we can interpet it as a relation on modes of $V$ with functional coefficients. 
In particular, all operators $T$ change vertex operator algebra elements by action either of 
$o(v)=v_{\wt v -1 }.$ or positive modes of $v(m).$, $m \ge 0$. 
Recall that for all $g$ the $n$-point functions 
possess certain modular properties with respect to groups depending on 
genus $g$ and modular space $\mathcal B^{(g)}$ parameters. 
Moreover, the reduction formulas \eqref{reduction}  
are used to prove modular invariance for $n$-point functions.   
Due to automorphic properties of $n$-point functions, 
 \eqref{conditions} can be also interpreted as relations among modular forms.  
The condition \eqref{reduction}  
defines a complex variety in ${\bf z}_n \in \C^{n}$ with non-commutative parameters 
${\bf v}_n \in V^{\otimes n}$.
As most identities 
(e.g., trisecant identity \cite{Fa, TZ1} and triple product identity \cite{K, MTZ, TZ2}) 
for $n$-point functions \eqref{conditions} has its algebraic-geometrical meaning. 
The condition \eqref{conditions} relates finite series of vertex operator algebra correlations functions on 
a genus $g$ Riemann surface 
with rational function coefficients (at genus $g=0$) \cite{FHL}, or elliptic functions (at genus $g=1$) 
\cite{Zhu, MT, MTZ}, or 
generalizations of classical elliptic functions (at genus $g \ge 2$) \cite{GT, TW}.  
Since in all cases $n$-point functions possess certain modular properties, we treat \eqref{conditions} as a source of 
new identities on modular forms at corresponding genus of Riemann surfaces. 
\section{Cohomology}
\label{cohomology}
In this section we compute the reduction cohomology defined by \eqref{buzova}--\eqref{pusto}. 
The main result of this paper is the following. 
\begin{proposition}
The $n$-th reduction cohomology of a vertex operator algebra $V$-module $W$ 
 is the space of analytical continuations of solutions  
$\F^{(g)}_{W, n}\left({\bf x}_n \right)$ 
to a Knizhnik-Zamolodchikov equation,    
 and provided by series of auxiliary functions
 recursively generated by reduction formulas \eqref{reduction} with $x_i \notin {\mathfrak V}_{i}$, 
for $1 \le i \le n$.    
\end{proposition}
\begin{remark}
 The first cohomology    
is given by 
 the space of transversal (i.e., with vanishing sum over $q$, $q'$)
 one-point connections $\F^{(g)}_{W, 1}\left(x_1 \right)$ provided by   
 coefficients in terms of series of generalizations of elliptic functions (cf. Appendix 6).   
The second cohomology is given by a space of   
generalized higher genus 
complex kernels corresponding to $V$ and $\Sigma^{(g)}$. 
\end{remark}
\begin{proof}
By definition \eqref{pusto}, 
the $n$-th reduction cohomology is defined by the subspace of $C^{(g)}_{n}(W)$ of functions 
$\F^{(g)}_{W, n}\left({\bf x}_n \right)$ satisfying 
\begin{eqnarray}
\label{poroserieroj}
&& \Big( \sum\limits_{l=1}^{l(g)}  
  f_1^{(g)}\left({\bf x}_{n+1}, l \right)  \;  T^{(g)}_l (v_{n+1})    
\nn
&& \qquad + 
 \sum\limits_{k=1}^{n} \sum\limits_{m \ge 0}  
  f_2^{(g)}\left( k, m \right) \;  T^{(g)}_k(v_{n+1}(m)) \Big) \;
 \F^{(g)}_{W, n}\left( {\bf x}_{n} \right)=0,  
\nn
\end{eqnarray} 
modulo the subspace of $C^{(g)}_{n}(W)$ $n$-point functions 
$\F_{W, n}^{(g)}\left({\bf x}'_n \right)$ resulting from:  
\begin{eqnarray}
\label{poroserieroj_2}
\F_{W, n}^{(g)} \left( {\bf x}'_n \right) &=& \Big( \sum\limits_{l=1}^{l(g)}   
  f_1^{(g)}\left({\bf x}'_n l \right)  \;  T^{(g)}_l   
\nn
&+& 
 \sum\limits_{k=1}^{n-1} \sum\limits_{m \ge 0}   
  f_2^{(g)}\left( k, m \right) \;  T^{(g)}_k(v'_n(m)) \Big) 
\; \F^{(g)}_{W, n-1}\left( {\bf x}'_{n-1} \right).  
\nn&&
\end{eqnarray}
We assume that, subject to other fixed parameters, $n$-point functions are completely 
determined by all choices ${\bf x}_n \in V^{\otimes n}\times \C^n$. 
Thus, the reduction cohomology can be treated as depending on set of ${\bf x}_n$ only
with appropriate action of endomorphisms generated by $x_{n+1}$. 
Consider a non-vanishing solution $\F^{(g)}_{W, n}\left({\bf x}_n  \right)$
 to \eqref{poroserieroj} for some ${\bf x}_n$. 
Let us use the reduction formulas \eqref{reduction} recursively for each $x_i$, $1 \le i \le n$ of ${\bf x}_n$ 
 in order to express 
$\F^{(g)}_{W, n}\left({\bf x}_n  \right)$ in terms of the partition function 
$\F^{(g)}_{W, 0}\left( B^{(g)}\right)$, , i.e., we obtain 
 \begin{equation}
\label{topaz}
\F^{(g)}_{W, n}\left({\bf x}_n  \right)= {\mathcal D}^{(g)}({\bf x}_n) \F^{(g)}_{W, 0}\left( B^{(g)}\right), 
\end{equation}
as in \cite{MT, MTZ, TZ, TZ1, TZ2}. 
It is clear that $x_i \notin {\mathfrak V}_{i}$ for $1 \le i \le n$, 
i.e., at each stage of the recursion procedure towards 
\eqref{topaz}, otherwise $\F^{(g)}_{W, n}\left({\bf x}_n  \right)$ would be zero. 
Thus, $\F^{(g)}_{W, n}\left({\bf x}_n  \right)$ is explicitly known and 
is repsented as a series of auxiliary functions ${\mathcal D}^{(g)}$  
depending on $V$, genus $g$, and moduli space parameters $B^{(g)}$. 
Consider now $\F^{(g)}_{W, n}\left({\bf x}'_n  \right)$ given by \eqref{poroserieroj_2}.  
It is either vanishes when $v_{n-i} \in {\mathfrak V}_{n-i}$, $ 2 \le i \le n$, or 
given by \eqref{topaz} with ${\bf x}'_n$ arguments. 

The way the reduction relations \eqref{reduction} were derived in \cite{MT, MTZ, GT, TW, BKT}   
is exactly the same as for the vertex operator algebra derivation \cite{KZ, TK} 
for the Knizhnik-Zamolodchikov equations. 
 The general idea is consider the double integration of 
 $\F^{(g)}_{W, n}\left({\bf x}_n \right)$ along small circles around two auxiliary variables 
with the action of reproduction kernels inserted. 
Then, these procedure leads to recursion formulas relating $\F^{(g)}_{W, n+1}$ and $\F^{(g)}_{W, n}$ with   
functional coefficients depending on the nature of the vertex operator algebra $V$, and the way 
a Riemann surface $\Sigma^{(g)}$ was constructed. 
Thus, in our context, \eqref{poroserieroj} is seen as a version of the Knizhnik-Zamolodchikov equation. 
In \cite{Y, MT, MTZ} formulas to $n$-point functions in various specific examples of $V$ and configuration 
of Riemann surfaces were explicitely and recursively obtained. 

In terms of $x_{n+1}$, by 
using \eqref{Ydefn}--\eqref{tornado}, we are able to transfer in \eqref{poroserieroj} the action 
of  
$v_{n+1}$-modes 
into an analytical continuation of 
$\F^{(g)}_{W, n}\left({\bf x}_n \right)$ 
multi-valued holomorphic functions to domains $D_{n} \subset  \Sigma^{(g)}$ with
 $z_{i} \neq z_{j}$ for 
$i\ne j$.
Namely, in \eqref{poroserieroj},
 the operators $T^{(g)}_l$ and $T^{(g)}_k(v_n(m))$ act by certain modes $v_{n+1}(m).$ of a vertex operator algebra 
element $v_{n+1}$ on ${\bf v}_n \in V^{\otimes n}$.  
Using vertex operator algebra associativity property \eqref{tornado} we express the action of 
of operators $T^{(g)}_l$ and $T^{(g)}_k(v_n(m))$ in terms of modes 
$v_{n+1}(m)$  inside vertex operators 
in  actions of $V$-modes on the whole vertex operator at expense of a shift of 
their formal parameters ${\bf z}_n$ by $z_{n+1}$, i.e., ${\bf z}'_n= {\bf z}_n + z_{n+1}$.   
Note that under such associativity transformations $v$-part of ${\bf x}_n$, i.e., ${\bf v}_n$ remains the
untouched. 
Thus, 
the $n$-th reduction cohomology of a $V$-module $W$ is given by the space of 
 analytical continuations of $n$-point functions 
$\F^{(g)}_{W, n}\left({\bf x}_n \right)$ with ${\bf x}_{n-1} \notin {\mathfrak V}_{n-1}$ 
that are solutions 
to the Knizhnik-Zamolodchikov equations \eqref{poroserieroj}. 
The above analytic extensions for the Knizhnik-Zamolodchikov equations generated by $x_{n+1}$
and 
with coefficients provided by functions 
$f_1^{(g)}\left({\bf x}_{n+1}, n \right)$ and $f_2^{(g)}\left( {\bf x}_{n+1}, k, m \right)$   
on the Riemann surface $\Sigma^{(g)}$)
\end{proof}
This result is coherent with considerations of \cite{TUY}.    
We illustrate the proposition above in cases of zero, one, two, and higher genera in section \ref{examples}.   
 One can make connection with the first cohomology of grading-restricted vertex operator algebras 
in terms of derivations, and  
 to the second cohomology in terms of square-zero extensions of $V$ by $W$ \cite{Huang2}. 

\medskip 
%
\noindent{\it Euler-Poincare formula. }
In \cite{Fuks}, for a Lie algebra, we find a celebrated formula relating sums over dimensions of 
chain complex spaces and dimensions of homologies for a graded infinite-dimensional Lie algebra. 
Suppose a $V$-module $W$ is endowed with a complex grading $W=\bigcup_{m\in \C} W_{(m)}$, i.e., one could consider 
$m= m' + \beta(m)$, where $m'\in \Z$, and $\beta(m)\in \C$ \cite{MTZ}.   
In our case all spaces $C^n_{(g)}(W)$ are infinite-dimensional, but, according to definition of 
a vertex operator algebra, $\dim W_{(m)} \le \infty$, $m \in \C$. Thus, for a fixed $g$, each $C^n_{(g)}(W)$ can be 
endowed with separation with respect to $m$.
Then the complex \eqref{buzova} decomposes into sum of complexes
\begin{equation}
0 \rmap  
C^{0}_{(g)}(m) \stackrel {\delta^{0}_{(g)}} {\rmap} C^1_{(g)}(m) 
\stackrel {\delta^{1}_{(g)}} {\rmap} 
\ldots \rmap C^{n-1}_{(g)}(m) \stackrel{\delta^{n-1}_{(g)}}{\rmap} C^n_{(g)}\rmap  \ldots. 
\end{equation}  
For $g \ge 0$, $l \ge 0$,  
let 
\[
q_{n, m}= \dim C_{(g)}^n
\left( V_{(m)} \right)= 
\dim 
\left\{
 \F^{(g)}_{W, n} 
\left( 
{\bf x}_n  
 \right): 
 {\bf x}_n \in V^{\otimes n}_{(m)}, m \in \C  \right\}. 
\]
%
Let 
\[
p_{m, n}=\dim H^n (V_{(m)})= 
\dim \left\{ \left({\rm Ker}\; \delta^n_{(g)} /{\rm Im} \; \delta^{n-1}_{(g)}\right)|_{V(m)},  
 {v}_{n+1}\in  V_{(m)} 
\right\},  
\]   
be dimensions of corresponding cohomology spaces. 
In our context we find that 
for fixed $m  \in \C$, and $n \ge 0$, 
\begin{equation}
\label{pofu}
\sum\limits_{n \ge 0 }^N (-1)^n \left( q_{n, m} - 
p_{n, m}\right)=0, 
\end{equation} 
Indeed, 
 let us consider relations of chain complex spaces and cohomology spaces for vertex operator algebra elements that 
belong to fixed $V_{(m)}$ grading subspace of $V$. 
Recall that according to the definition of a vertex operator algebra (see Appendix \ref{vertex})  
subspaces $W_{(m)}$ are finite dimensional.  
Thus, as a functional space,
 $\left\{ \F^{(g)}_{W, n}\left({\bf x}_n \right), {\bf v}_n \in V^{\otimes n}_{(m)}\right\}$ 
is finite-dimensional. 
Consider now the cohomology spaces $\left\{ H_{(g)}^n, {\bf v}_{n+1} \in V^{\otimes (n+1)}_{(m)}\right\}$. 
Let us fix $v_{n+1} \in V_{(m)}$. 
For the same reason as above, ${\rm rank} \left( {\rm Ker} \; \delta^n_{(g)} \right) < \infty$, 
${\rm rank} \left( {\rm Im} \; \delta^{n-1}_{(g)} \right)  < \infty$,
as subspaces of $n+1$- and $n$-point functions. 
Thus, ${\rm rank} H^n_{(g)}|_{V_{(m)}} < \infty$. 
Using then the standard Eulre-Poincare formula \cite{Fuks} considerations for rank we obtain \eqref{pofu}.  
The formula \eqref{pofu} has deep number-theoretical meanings as equality of generating functions 
for series expansions for automorphic forms. 
\section{Examples}
\label{examples}
The reduction cohomology depends on actual coefficients of the Knizhnik-Zamolodchikov
 equations \eqref{poroserieroj}.   
Note that for $n=0$, $\F^{(g)}_{W, 0}\left(B^{(g)}\right)$, $g\ge 0$, 
is the partition function (or graded dimension for $g=1$),  
for a vertex operator algebra $V$ module $W$. 
In this section we provide examples of vertex operator algebras considered 
on Riemann surfaces of various genus $g$. 
\subsection{Vertex operator algebra $n$-point functions on the sphere} 
\label{sphere}
For $u$, ${\bf v}_{n}\in V$, and a homogeneous $u \in V$, we find the formula  
\cite{FHL, FLM}  
for the $n$-point functions  
\[
\F^{(0)}_{W, n}({\bf x}_{n})= 
 \langle u^{\prime }, {\bf  Y ( x}_n{\bf )}  
u \rangle, 
\]
where 
\[
{\bf  Y} ( {\bf x}_n ) =Y(x_1) \ldots Y(x_n),  
\]
The partition function is given by 
\[
\F^{(0)}_{W, 0}=\langle u'_{(a)}, u_{(b)} \rangle= \delta_{a,b}.   
\]
The reduction operators of \eqref{poros} are \cite{FHL} 
\begin{eqnarray}
\label{zhu_reduction_genus_zero_1} 
 H_1^{(0)}({\bf x}_{n+1})\; \F^{(0)}_{W, n}({\bf x}_{n}) &=& T_1(o(v)) \; \F^{(0)}_{W, n}({\bf x}_{n}),   
\\
\nonumber
 H_2^{(0)}({\bf x}_{n+1})\; \F^{(0)}_{W, n}({\bf x}_n) 
&=&\sum\limits_{k=1}^{n} \sum\limits_{m \ge 0}    
 f_{wt(v_{n+1}), m}(z_{n+1}, z_r) \; T_k( v(m) ) \; \F^{(0)}_{W, n}({\bf x}_{n}), 
\end{eqnarray}
where 
 we define $f^{(0)}_{wt(v, m} (z, w)$ is a rational function defined by 
\[
f^{(0)}_{n,m}(z,w) 
=
\frac{z^{-n}}{m!}\left(\frac{d}{dw}\right)^m \frac{w^n}{z-w}, 
\] 
\[
\iota_{z,w}f^{(0)}_{n,m}(z,w) 
=
\sum\limits_{j\in {\mathbb N}}\left( { n+j \atop m}\right) z^{-n-j-1}w^{n+j-1}. 
\]
Let us take $x_{n+1}$ as the variable of expansion.  
The $n$-th reduction cohomology $H^n_{(0)}(W)$ is given the space of solutions to the Knizhnik-Zamolodchikov 
equation \eqref{perdolo1}  
 of the $n$-point function $\F^{(0)}_{W, 0}({\bf x}_n)$ 
(not given by $\delta^{n-1}_{0)} \F^{(0)}_{W, n-1}({\bf x}_{n-1})$)  
\begin{eqnarray}
{\bf Y}({\bf x}_n) &\mapsto& o(v_{n+1})\; {\bf Y}({\bf x}_n),  
\nn
 {\bf x}_k & \mapsto& v_{(m)}. {\bf x}_k, 
\end{eqnarray}
(i.e., generated by $x_{n+1}$-endomorphisms)  
of solututions $\F^{(0)}_{V, n}$ of the Knizhnik-Zamolodchikov equation 
\begin{eqnarray}
\label{perdolo}
 0 =  \F_{W, n}^{(0)}\left(o(v_{n+1}) \; {\bf x}_{n}\right) 
 +\sum\limits_{k=1}^{n} \sum\limits_{m \ge 0}  
 f^{(0)}_{wt(v_{n+1}), m} (z_{n+1}, z_k) \;  \F^{(0)}_{W, n}((v_{n+1}(m))_k.{\bf x}_n),     
\end{eqnarray}
 with rational function coefficiens $f^{(0)}_{wt(v_{n+1}), m, k}(z_{n+1}, z_k)$,  
modulo the space of $n$-point functions obtained by the recursion procedure. 
Using \eqref{Ydefn}--\eqref{tornado} we obtain from \eqref{perdolo}, and the standard Virasoro algebra 
representation $L_V(0)=z_{n+1} \partial_{z_{n+1}}$, we obtain  
\begin{eqnarray}
\label{perdolo1}
 \left(  \partial_{z_{n+1}}   
 + \sum\limits_{k=1}^{n} 
   \widetilde{f}^{(0)}_{wt(v_{n+1}), m} (z_{n+1}, z_k) \right) \;  \F^{(0)}_{W, n}({\bf x}_n + (z_{n+1})_k)=0,     
\end{eqnarray}
which is the Knizhnik-Zamolodchikov equation on analytical continuation of  
$\F^{(0)}_{W, n}({\bf x}_n + (z_{n+1})_k)$ with a different $\widetilde{f}^{(0)}$.  
Using the reduction formulas \eqref{reduction} we obtain 
\[
\F^{(0)}_{W, n}({\bf x}_n + (z_{n+1})_k)= {\mathcal D}^{(0)}({\bf x}_{n+1}), 
\]
where ${\mathcal D}^{(0)}({\bf x}_{n+1})$ is given by the series of 
rational-valued functions in ${\bf x}_{n+1} \notin {\mathfrak V}_n$ resulting from the recursive procedure 
starting from $n$-point function to the partition function.  
Thus, in this example, the $n$-th cohomology 
is the space of analytic extensions \cite{FHL, FLM} of rational function solutions to 
 the equation \eqref{perdolo} with rational function coefficients.   
\subsection{Vertex operator algebra $n$-point functions on the torus}
\label{torus}
In order to consider modular-invariance of $n$-point functions at genus one, Zhu
introduced \cite{Zhu} a second square-bracket VOA 
$(V,Y[.,.],\mathbf{1}_V, \tilde{\omega})$ 
associated to a given VOA $(V,Y(.,.),\mathbf{1}_V,\omega )$.
 The new square bracket vertex operators are 
\[
Y[v,z]=\sum_{n\in \mathbb{Z}}v[n]z^{-n-1}=Y(q_{z}^{L(0)}v,q_{z}-1),
\]
with $q_{z}=e^z$, 
while the new conformal vector is 
\[
\tilde{\omega}
=\omega -\frac{c_V}{24}\mathbf{1}_V.
\]
For $v$ of $L(0)$ weight $wt(v)\in \mathbb{R}$ and $m\geq 0$, 
\[
v[m] 
=
m!\sum\limits_{i\geq m}c(wt(v),i,m)v(i),
\]
where 
\[
\sum\limits_{m=0}^{i}c(wt(v),i,m)x^{m} 
=
  \left( wt(v)-1+x \atop i \right).
\]
For ${\bf v}_{n}\in V^{\otimes n}$ the genus one $n$-point function \cite{Zhu} has the form  
\[
F^{(1)}_n({\bf x}_n; \tau)
=
Tr_V\left({\bf Y} \left({\bf q^{L(0)} v}_n, {\bf q}_n\right) \; q^{L(0)-c_V/24}\right), 
\]
for $q=e^{2\pi i \tau}$ and $q_{i}=e^{z_{i}}$, where $\tau$ is the torus modular parameter, 
and $c_V$ is the central charge of $V$-Virasoro algebra.    
Then the genus one Zhu recursion formula is given by the following \cite{Zhu}. 
For any $v_{n+1}$, ${\bf v_{n}}\in V^{\otimes n}$ we find for an $n+1$-point function 
\begin{eqnarray}
 H_1^{(1)}({\bf x}_{n+1})  &=& F^{(1)}_{n} \left( o(v_{n+1}) \; {\bf x}_{n} ; \tau \right)
\\
  H_2^{(1)}({\bf x}_{n+1}) &=&\sum\limits_{k=1}^{n}\sum\limits_{m \geq 0} 
 P_{m+1}
(z_{n+1}-z_{k},\tau )
F^{(1)}_V ( (v[m])_k. \; {\bf x}_n; \tau).  
\nonumber
\label{zhu_reduction_genus_one}
\end{eqnarray}
In this theorem $P_{m}(z,\tau)$ denote higher Weierstrass functions defined by 
\[
P_{m}
 (z,\tau )=\frac{(-1)^{m}}{(m-1)!}\sum\limits_{n\in {\mathbb Z}_{\neq 0}
 } 
\frac{n^{m-1}q_{z}^{n}}{1-q^{n}}.
\]  
The partition function  
\begin{eqnarray}
\label{nova2}
\F^{(1)}_{W, 0} ( \tau) &=&  {\rm Tr}_V\left( q^{L_V(0)-c/24}\right), 
\end{eqnarray}
 is called the graded dimension for $V$.
Let us fix ${\bf x}_{n+1} \in V^{\otimes n}\times \C^n$.
The $n$-th reduction cohomology 
$H^n$ is the space of extensions of the Knizhnik-Zamolodchikov equation solutions 
\begin{eqnarray}
\label{popaz} 
0 = F^{(1)}_V \left( o(v_{n+1})\;  {\bf x_{n}}; \tau \right)
  + \sum\limits_{k=1}^{n+1} \sum\limits_{m\geq 0} 
 P_{m+1}
(z_{n+1}-z_{k},  \tau ) 
F^{(1)}_{n} ((v_{n+1}[m])_k. {\bf x}_n; \tau).    
\end{eqnarray}
This space can be described by complex kernels (i.e., prime forms \cite{Mu}) given by 
sums of elliptic functions.  
\subsection{Genus two $n$-point functions}   
\label{corfu}
In this subsection we recall \cite{GT} the construction and reduction formulas for
  vertex operator algebra 
$n$-point functions on a genus two Riemann surface formed in the sewing procedure of two torai
 due to \cite{Y}. 
Let us 
assume that a vertex operator algebra $V$ is of strong-type. Thus,  it possess  
a non-degenerate 
bilinear form.  
For a $V$--basis $\{u^{(a)}\}$ we define the dual basis $\{\overline{u}^{(a)}\}$ with respect 
to the bilinear form  
where
\begin{align*}
\langle \overline{u}^{(a)}, u^{(b)} \rangle_{\mathrm{sq}} = \delta_{ab}.
\end{align*}
\begin{definition}
The
 {genus two partition function} (zero-point function) for $V$ is defined by
\begin{align}
\label{eq:Zdef}
\F^{(2)}_{V, 0}\left(B^{(2)}\right) = \sum_{r \geq 0} \epsilon^r \sum_{\overline{u}\in V_{[r]}} 
\F^{(1)}_{V, 1} (u; \tau_1) \; \F^{(1)}_{V, 1}({u}; \tau_2),
\end{align}
where $B^{(g)}=(\tau_1, \tau_2, \epsilon)$ in the $\epsilon$-sewing procedure 
for constructing a genust two Riemann surface \cite{Y, MT, TZ1},  
 and 
the internal sum is taken over any $V_{[n]}$--basis, and $\overline{u}$
 is the dual of $u$ with respect to the bilinear form on $V$. 
\end{definition}
We then recall \cite{GT}     
a formal genus two reduction formulas for all
 $n$-point functions.  
\begin{definition}
Let  $v_{n+1} \in  V$ be inserted at $x_{n+1} \in \widehat{\Sigma}^{(1)}_1$, ${\bf v}_k \in V^{\otimes k}$
 be inserted at   
${\bf y}_k \in \widehat{\Sigma}^{(1)}_1$ and  
${\bf v}'_l \in V^{\otimes l}$  
be inserted at ${\bf y}'_l \in\widehat{\Sigma}^{(1)}_2$,  
a punctured torus 
\begin{align}\label{Sahat}
\widehat{\Sigma}^{(1)}_{a}=
\Sigma^{(1)}_{a}\backslash\left \{z_{a},\left\vert z_{a}\right\vert \leq |\epsilon |/r_{\bar{a}}\right\},
\end{align}
where we here  
we use the convention  
\begin{align}\label{abar}
\overline{1}=2,\quad \overline{2}=1.
\end{align}
We consider  the corresponding genus two $n$-point function 
\begin{eqnarray}
&&\F^{(2)}_{V, n} \left(x_{n+1}; {\bf x}_k;  
{\bf x}'_l;  B^{(2)} 
\right)  
\nn
&& \qquad  = 
\sum_{r\geq 0} \epsilon^r \sum_{u\in V_{[r]}} 
\F^{(1)}_{V, k+1} \left(Y[x_{n+1}] \;  
{\bf Y[x_k]}_k\; 
u; \tau_1 \right) \; 
\F^{(1)}_{V, l} \left( 
 {\bf Y[x'_l]}_l \;  
\overline{u}; \tau_2 \right).  
\label{Znvleft1}
\end{eqnarray}
where the sum as in \eqref{eq:Zdef}.
\end{definition}
\begin{remark}
\label{rem:basis}
\eqref{eq:Zdef} is independent of the choice of $V$--basis. 
One could also write a 
 similar expression by inserting the $Y[x]$ vertex operator  at
$z \in\widehat{\Sigma}^{(1)}_2$ on the right hand 
side of \eqref{Znvleft1}.  
\end{remark}
First, one defines the functions $\F^{(2)}_{V, n, a}$ for $a\in \{1,2\}$,  by
\begin{eqnarray*}
\label{eq:Oadef}
 \F^{(2)}_{V, n, 1} \left(  {\bf x}_{n+1}; B^{(2)} \right) 
  &=&\sum_{r\geq 0} \epsilon^r \sum_{u\in V_{[r]}} 
{\rm Tr}_V \left( o(v_{n+1}) \;  
{\bf Y( q_{z_k}^{L(0)}   } {\bf v}_k, {\bf q_{z_k})} \; 
q_1^{L_V(0) - c_V/24} \right) 
\nn
  && \qquad \cdot \F^{(1)}_{V, n-k} 
\left(
{\bf Y} [    {\bf x}_{k+1, n}  ]  \; 
\overline{u}; \tau_2
\right), 
\nn
&&
\end{eqnarray*}
\begin{eqnarray*}
  \F^{(2)}_{V, n, 2} \left({\bf x}_{n+1}; B^{(2)}   
\right) 
&=&
\sum_{r\geq0} \epsilon^r \sum_{u\in V_{[r]}}
\F^{(1)}_{V, k}
\left(
{\bf Y[x_k]} \;
u;  \tau_1\right) \; 
\nn
&&
\; 
\cdot {\rm Tr}_V \left( o(v_{n+1})\;  
{\bf Y(q_{z_{k+1, n}}^{L(0)}} {\bf v_{k+1}}, {\bf q_{z_{k+1}}) }\; 
 q_2^{L_V(0) - c_V/24} \right),  
\end{eqnarray*}
and 
$\F^{(2)}_{V, n, 3} \left(   {\bf x}_{n+1}; B^{(2)} \right)= \mathbb{X}_1^\Pi$ of \eqref{eq:Xadef}, 
We also define 
\begin{definition}
\label{calFforms}
Let $f^{(2)}_{a}(p; z_{n+1})$, for $p \ge 1$, and $a=1$, $2$  be given by 
\begin{equation}
f^{(2)}_{a} (p;z_{n+1})=   
1^{\delta_{ba}} + (-1)^{p\delta_{b\abar} } \epsilon^{1/2} \left( 
\mathbb{Q}(p; z_{n+1}) \left( \widetilde{\Lambda}_{\overline{a}} \right)^{\delta_{ba} } \right) (1), 
\label{eq:calFadef} 
\end{equation}
for $z_{n+1}\in\widehat{\Sigma}^{(1)}_b$. 
Let
$f^{(2)}_3 (p; z_{n+1})$, for $z_{n+1} \in {\Sigma}^{(1)}_a$,  be an infinite row vector  
given by
\begin{equation}
f^{(2)}_3(p; z_{n+1})=    
\left( \mathbb{R}(z_{n+1}) + 
\mathbb{Q}(p;z_{n+1}) \left(\widetilde{\Lambda}_{\overline{a}}\Lambda_{a} 
+ \Lambda_{\overline{a}}\Gamma \right) \right)\Pi.
\label{eq:calFPidef}
\end{equation}
\end{definition}
In \cite{GT} it is proven that 
 the genus two $n=k+l$-point  function for a quasi-primary vector $v_{n+1}$ of weight 
$p=\wt{[v_{n+1}]}$ inserted at 
$x \in\widehat{\Sigma}^{(1)}_1$, and general vectors ${\bf v_k}$ and ${\bf v}'_l$  inserted at 
${\bf x}_k\in \widehat{\Sigma}^{(1)}_1$,  
respectively, has the following operators $H_1^{(2)}({\bf x}_{n+1})$ and $H_2^{(2)}({\bf x}_{n+1})$: 
\begin{eqnarray}
&& H^{(2)}_1 ({\ x}_{n+1} ) \; \F^{(2)}_{V, n} \left({\bf x}_n; B^{(2)}
\right)  
 = \sum\limits_{l=1}^3  
f^{(2)}_{l}(p; z_{n+1}) \; \F^{(2)}_{V, n, l} \left({\bf x}_{n+1}; B^{(2)} \right),   
\nn
&& 
H^{(2)}_2({x}_{n+1}) \; \F^{(2)}_{V, n} \left( {\bf x}_n; B^{(2)}\right)= \sum_{i=1}^n \sum_{j\geq0} 
\mathcal{P}_{j+1}(p; z_{n+1}, z_i) \;  
\F^{(2)}_{V, n} \left( (v_{n+1}[j])_i. {\bf x}_n;  B^{(2)} \right),   
\nn
&&
\label{eq:npZhured}
\end{eqnarray}
with $\mathcal{P}_{j+1}(p;x,y)$ of \eqref{eq:PN21j}. 
A similar expression corresponds to $v_{n+1}$ inserted on $x\in \widehat{\Sigma}^{(1)}_2$.
\begin{remark}
The functions 
$\mathcal{P}_{j+1}(p;x,y)$ are holomorphic for $x \neq y$ on the sewing domain 
 in the cases $p=1$, $2$.  
\end{remark}
\subsection{Vertex operator algebra $n$-point functions and reduction formulas in genus $g$ Schottky parameterization} 
\label{corfug}
In this section we recall \cite{TZ, TW, T2} the construction and reduction relations 
 for vertex operator algebra $n$-point functions  
defined on a genus $g$ Riemann surface formed in the Schottky parameterization. 
In particular, 
the formal partition and $n$-point correlation functions for a vertex operator 
algebra associated to a genus $g$ Riemann surface $\Sigma^{(g)}$ 
 are introduced in the Schottky 
 scheme.  
All expressions here are functions of formal variables $w_{\pm a}$, $\rho_{a} \in \C$ 
 and vertex operator parameters. 
Then we recall   
 the genus $g$ reduction formula with universal coefficients that have a geometrical meaning and 
are meromorphic on a Riemann surface 
$\Sigma^{(g)}$. 
 These coefficients are  
generalizations of the 
elliptic Weierstrass functions \cite{L}. 
For a $2g$ vertex operator algebra $V$ states 
\[
\bm{b}=(b_{-1}, b_{1}; \ldots;  b_{-g} ;b_{g}),  
\]
and corresponding local coordinates 
\[
\bm{w}= (w_{-1}, w_{1}; \ldots; w_{-g},w_{g}),  
\]
of $2g$ points $(p_{-1}, p_1; \ldots;  p_{-g}, p_g)$ on the Riemann sphere $\Sigma^{(0)}$,  
consider the genus zero $2g$-point correlation function
\begin{align*}
Z^{(0)}(\bm{b,w})=&Z^{(0)}(b_{-1},w_{-1};b_{1},w_{1};\ldots;b_{-g},w_{-g};b_{g},w_{g})
\\
=&\prod_{a\in\Ip}\rho_{a}^{\wt(b_{a})}Z^{(0)}(\bbar_{1},w_{-1};b_{1},w_{1};\ldots;\bbar_{g},w_{-g};b_{g},w_{g}).
\end{align*}
where $\Ip=\{1,2,\ldots,g\}$. 
Let 
\[
\bm{b}_{+}=(b_{1},\ldots,b_{g}), 
\]
 denote an element of a $V$-tensor product $V^{\otimes g}$-basis with dual basis 
\[
 \bm{b}_{-}=(b_{-1}, \ldots,b_{-g}), 
\]
 with respect to the bilinear form 
$\langle \cdot, \cdot\rangle_{\rho_{a}}$ (cf. Appendix \ref{vertex}). 

 Let $w_{a}$ for $a\in\I$ be $2g$ 
 formal variables. One identify them  
 with the canonical Schottky parameters (for detailes of the Schottky construction, see \cite{TW, T2}).    
 One can define the genus $g$ partition function as
\begin{align}
\label{GenusgPartition}
Z_{V}^{(g)} =Z_{V}^{(g)}(\bm{w,\rho})
=\sum_{\bm{b}_{+}}Z^{(0)}(\bm{b,w}),
\end{align}
for 
\[
(\bm{w,\rho})=(w_{\pm 1},  \rho_{1}; \ldots; w_{\pm g},\rho_{g}). 
\]
Now  
we recall 
the formal reduction formulas for all
 genus $g$  Schottky $n$-point functions.
%
One 
defines the genus $g$  formal $n$-point function for $n$ vectors 
${\bf v}_{n} \in V^{\otimes n}$ inserted at ${\bf y}_{n}$ by 
\begin{align}\label{GenusgnPoint}
Z_{V}^{(g)}(\bm{v,y}) =Z_{V}^{(g)}(\bm{v,y};\bm{w,\rho})
=
\sum_{\bm{b}_{+}}Z^{(0)}(\bm{v,y};\bm{b,w}),
\end{align}
where 
\[
Z^{(0)}(\bm{v,y};\bm{b,w})=Z^{(0)}(v_{1},y_{1};\ldots;v_{n},y_{n};b_{-1},w_{-1};\ldots;b_{g},w_{g}).
\] 
Let $U$ be a vertex operator subalgebra of $V$ where $V$ has a $U$-module decomposition
\[
V=\bigoplus_{\alpha\in A}
{W}_{\alpha},
\]
  for  $U$-modules $
{W}_{\alpha}$ and some indexing set $A$. 
Let 
\[
{W}_{\bm{\alpha}}=\bigotimes_{a=1}^{g} 
{W}_{\alpha_{a}}, 
\]
 denote a tensor product of $g$ modules 
\begin{align}\label{eq:Z_Walpha}
Z_{
{W}_{\bm{\alpha}}}^{(g)}(\bm{v,y}) =\sum _{\bm{b_{+}}\in 
{W}_{\bm{\alpha}}} Z^{(0)}(\bm{v,y};\bm{b,w}),
\end{align}
where here the sum is over a basis $\{\bm{b}_{+}\}$ for  $ 
{W}_{\bm{\alpha}}$. 
It follows that
\begin{align}
\label{eq:Z_WalphaSum}
Z_{V}^{(g)}(\bm{v,y})=\sum_{\bm{\alpha}\in\bm{A}}Z_{ 
{W}_{\bm{\alpha}}}^{(g)}(\bm{v,y}),
\end{align}
where the sum ranges over $\bm{\alpha}=(\alpha_{1},\ldots ,\alpha_{g}) \in \bm{A}$,  for $\bm{A}=A^{\otimes{g}}$. 
Finally, one defines corresponding formal $n$-point correlation differential forms 
\begin{align}
\label{tupoy}
\F_{V}^{(g)}(\bm{v,y}) &=Z^{(g)}(\bm{v,y}) \; \bm{dy^{\wt(v)}},
\nn
\F_{
{W}_{\bm{\alpha}}}^{(g)}(\bm{v,y}) &=Z_{
{W}_{\bm{\alpha}}}^{(g)}(\bm{v,y})\; \bm{dy^{\wt(v)}},
\end{align}
where  
\[
\bm{dy^{\wt(v)}}=\prod_{k=1}^{n} dy_{k}^{\wt(v_{k})}.
\]
Recall notations and 
identifications given in Appendix \ref{derivation}. 
In \cite{TW} they prove that 
 the genus $g$ $(n+1)$-point formal differential 
$\F_{W_{\bm{\alpha}}}^{(g)}(x; \bm{v,y})$,  for $x_{n+1}= (v_{n+1}, y_{n+1})$, 
for a quasiprimary vector $v_{n+1} \in U$ of weight $\wt(v_{n+1})=p$ inserted at a 
 point
$p_0$, with the coordinate $y_{n+1}$,   
 and general vectors ${\bf v_{n}}$ inserted at points 
 ${\bf p}_n$ with coordinates 
${\bf y}_{n}$ satisfies the recursive identity for ${\bf x}_n= (\bm{v, y})$
\begin{eqnarray}
\label{eq:ZhuGenusg}
\F_{{
{W, n+1}_{\bm{\alpha}}}}^{(g)}\left(x_{n+1};  
{\bf x}_n 
\right) & = & \left(H_1^{(g)} + H_2^{(g)}\right) \; \F^{(g)}_{W, n}\left({\bf x}_n \right),   
\\
H_1^{(g)}({\bf x}_{n+1}) \F^{(g)}_{W, n} \left({\bf x}_n \right) &=& \sum_{a=1}^{g} \Theta_{a}(y_{n+1}) 
\; O^{W_{\bm{\alpha}}}_{a}\left(v_{n+1}; {\bf x}_n\right),   
\nn
H_2^{(g)}({\bf x}_{n+1})\F^{(g)}_{W, n+1} \left({x_{n+1}}  \right) &=&  
\sum_{k=1}^{n}\sum_{j\ge 0}\del^{(0,j)} \; \Psi_{p}(y_{n+1},y_{k})\; 
\nn
&& \qquad 
\cdot \F_{{
{W, n}_{\bm{\alpha}}}}^{(g)} 
\left( (u(j))_k. {\bf x_n} \right)\; dy_{k}^{j}\; 
\nonumber
\end{eqnarray}
Here $\del^{(0,j)}$ is given by     
\[
\del^{(i,j)}f(x,y)=\del_{x}^{(i)}\del_{y}^{(j)}f(x,y), 
\]
for a function $f(x,y)$, and $\del^{(0,j)}$ denotes partial derivatives with respect to $x$ and $y_j$. 
The forms 
$\Psi_{p}(y_{n+1}, y_{k} )\; dy_{k}^{j}$ given by 
 \eqref{psih}, 
 $\Theta_{a}(x)$ is of \eqref{thetanew}, and 
 $O^{
{W}_{\bm{\alpha}}}_{a}(v_{n+1}; {\bf x}_n)$ of \eqref{oat}.  
Similar definitions as in previous subsections can be formulated for a vertex operator algebra module. 
\subsection{Vertex operator cluster algebras}
\label{cluster}
Finally, we would like to describe the example of vertex operator cluster algebras introduced in \cite{Zu} 
and reveal their cohomological nature. 
The condition \eqref{conditions} for mutation \eqref{hruza} to be involutive  
gives us a cohomological condition 
(vanishing square of the $n$-th coboundary operator) 
on modular forms.  
Cluster algebras introduced in\cite{FZ1} 
have 
numerous applications in various areas of mathematics 
\cite{FG1, FG2, FG3, FG4, GSV1, GSV2, GSV3, FST, KS, N1, N2, DFK, GLS, HL, Ke1, Ke2, Na, Sch, N1, N2, Sch}.  
In this subsection we recall \cite{Zu} definition of a vertex operator cluster algebra 
and show its cohomological nature. 
 Let us fix a vertex operator operator algebra $V$. 
Chose $n$-marked points $p_i$, $i=1, \ldots, n$ on a genus g compact 
Riemann surface formed by one of procedures \cite{Y, TZ} mentioned in Introduction.  
In the vicinity of each marked point $p_i$ define a local coordinate $z_i$ 
with zero at $p_i$.  
For $n$-tuples of elements ${\bf x}_n$ 
let us denote 
\[
{\mathcal Y}({\bf x}_n)=(Y(x_1), \ldots, Y(x_n)).  
\] 
\begin{definition}
\label{seeda}
We define a vertex operator cluster algebra seed   
\begin{equation}
\label{vertex_cluster}
\left( {\mathbf v}_n, {\mathcal Y}({\bf x}_n), \F^{(g)}_{W, n} \left({\mathbf x}_n  \right) \right).   
\end{equation}
\end{definition}
The mutation is defined as follows: 
\begin{definition}
\label{muta}
For ${\mathbf v}_n$, we define the mutation ${\mathbf v}^\prime_n$ of 
 ${\mathbf v}_n$ in the direction $k \in 1, \ldots, n$ 
as 
\begin{equation}
\label{v-mutation}
{\mathbf v}^\prime_n= \mu_k( x_{n+1}, m) {\mathbf v}_n=  
\left((F_k(v_{n+1}(m)))_{k} \; {\bf v}_n \right), 
\end{equation}
for some $m \ge 0$, and $V$-valued functions $F^{(g)}_k(v_{n+1}(m))$, depending on genus $g$ of the Riemann surface. 
Note that due to the property (\ref{lowertrun})  
we get a finite number of terms as a result of the action of $v_{n+1}(m)$ on $v_k$, $1 \le k \le n$. 
For the $n$-tuple of vertex operators we define
\begin{eqnarray}
\label{z-mutation}
{\mathcal Y} \left({\bf x}^\prime_n \right) 
= \mu_k(x_{n+1}, m) \;{\mathcal Y} \left({\bf x}_n \right)  
= \left( 
Y \left( G^{(g)}_k(v_{n+1}(m)). {\bf x}_n\right) \right),     
\end{eqnarray}
where $G^{(g)}_k(v_{n+1}(m))$ are other $V$-valued functions.  
For $u \in V$, $w\in \C$,   
 the mutation $\mu(x_{n+1})$ of $\F^{(g)}_n \left({\mathbf x}_n \right)$,   
\begin{eqnarray}
 \F_{W, n}^{(g)} {}' \left(   {\mathbf x}^\prime_n \right) 
 = \mu
(x_{n+1})\; \F^{(g)}_{W, n} \left({\mathbf x}_n \right),  
\label{matrix_rule}
\end{eqnarray}
where $\mu(x_{n+1})$ is given by the coboundary operator 
\[
\mu(x_{n+1}) =\delta^n_{(g)}, 
\] 
defined by summation  
over mutations in all possible directions $k$,  $1 \le k \le n$, with   
auxiliary functions
$f_1^{(g)}\left(x_{n+1}, l \right)$ and     
  $f_2^{(g)}\left( m, k  \right)$, $1 \le l \le l(g)$, $1 \le k \le n$, 
$m \ge 1$, 
\begin{eqnarray}
\label{matrix_rule_n_plis_one}
 &&
 \F^{(g)}_{W, n} {}' \left({\mathbf x}^\prime \right)  
%
= \delta^n_{(g)} 
\F^{(g)}_{W, n} \left({\bf x}_n \right).  
\end{eqnarray}
We also require
the involutivity condition
\begin{equation}
\label{hruza}
\mu(x_{n+1}) \; \mu(x_{n+1}) = {\rm Id},  
\end{equation}
define the mutation of the seed (\ref{vertex_cluster}) in the direction $k \in 1, \ldots, n$   
of $x_{n+1} \in V$.  
We denote by $\mathcal C^n_{(G)}$ the space of seeds \eqref{vertex_cluster} for particular 
$n \ge 1$, 
satisfying condition \eqref{hruza}. 
\end{definition}
The involutivity condition can be expressed in terms of coboundary operator $\delta^n_{(g)}$ 
as \eqref{conditions} with $v_{n+2} = v_{n+1}$ with extra term ${\rm Id}$ on the right hand side.  
\begin{definition}
For a fixed $g$, 
 definitions \ref{seeda}--\ref{muta} and involutivity   
condition \eqref{hruza} for mutation   
determine the structure of a vertex operator cluster algebra ${\mathcal C}{\mathcal G}_n$   
of dimension $n$.  
We call the full vertex operator cluster algebra the union 
$\bigcup_{n \ge 0} \; {\mathcal C}{\mathcal G}_n$.  
\end{definition} 
It is naturally graded by $n$. 
 of ${\mathcal C}{\mathcal G}_n$.  
The cohomology of vertex operator cluster algebras will be considered elsewhere. 
In \cite{Zu} we have proven the following 
\begin{proposition}
\label{pizda}
For a vertex operator algebra $V$ such that $\dim V_k=1$, $k\in \Z$, and 
with $u=\one_V$,  $w\in \C$, and 
\[
F_k(u(m)).v=G_k(u(m)).v= \xi_{u, v} u[-1].v, 
\]
\[
  T_k^{(g)}(u(m))= u[m],  
\]
for $ m \ge 0$, and $\xi_{u, v}\in \C$, $\xi_{u, v}^2=1$, depending on $u$ and $v$, 
in \eqref{v-mutation}, \eqref{z-mutation}, and \eqref{matrix_rule_n_plis_one},  
the mutation 
\[
\mu 
= 
\left( \mu_k(\one_V, -1), \mu_k(\one_V, -1),  \mu(\one_V, w) \right), 
\]
\begin{equation}
\label{mutation}
\left({\mathbf v}^\prime_n, {\mathcal Y} ({\bf x}^\prime_n),  
\F^{(g)}_{W, n}{}'({\mathbf x}^\prime_n )\right)   
=\mu 
\;\left( {\mathbf v}, {\mathcal Y}({\bf x}_n), \F^{(g)}_{W, n} \left({\mathbf x}_n  \right) \right),  
 \end{equation}
defined by (\ref{v-mutation}), (\ref{z-mutation}), (\ref{matrix_rule_n_plis_one}) 
is an involution,   
i.e., 
\[
\mu
\; \mu 
= {\rm Id}. \;  
\]
\hfill $\square$
\end{proposition}
%
\section{Appendix: Vertex operator algebras} 
\label{vertex}
In this subsection we recall the notion of a vertex operator algebra \cite{B, DL, FHL, FLM, K, LL, MN}.  
\begin{definition}
A vertex operator algebra is determined by a quadruple  
$(V,Y,\mathbf{1}_V,\omega)$,  
where 
is a 
linear space 
endowed with a $\mathbb Z$-grading with 
\[
V=\bigoplus_{r \in {\mathbb Z}} V_{r},   
\]
with $\dim V_{r}<\infty$.
  The state  ${\mathbf 1}_V \in V_{0}$, $\mathbf{1}_V \not= 0$,  is the vacuum vector and 
$\omega\in V_{2}$ is the conformal vector with properties described below.
The vertex operator $Y$ is a linear map 
\[
Y: V\rightarrow \mathrm{End}(V)[[z,z^{-1}]], 
\]
for formal variable $z$ so that for any vector $u\in V$ we have a vertex operator  
\begin{equation}
Y(u,z)=\sum_{n\in {\mathbb Z}}u(n)z^{-n-1}.  
\label{Ydefn}
\end{equation}
The linear operators (modes) $u(n):V\rightarrow V$ satisfy creativity 
\begin{equation}
Y(u,z){\mathbf 1}_V= u +O(z),
\label{create}
\end{equation}
and lower truncation 
\begin{equation}
u(n)v=0,
\label{lowertrun}
\end{equation}
conditions for each $u$, $v\in V$ and $ n\gg 0$. 
Finally, the vertex operators satisfy the Jacobi identity
\begin{eqnarray}
\nonumber
 &&z_0^{-1}\delta\left( \frac{z_1 - z_2}{z_0}\right) Y (u, z_1 )Y(v , z_2)    
  - 
z_0^{-1} \delta \left( \frac{z_2 - z_1}{-z_0}\right) Y(v, z_2) Y(u , z_1 ) 
\\&&
\label{vertex operator algebraJac}
 \qquad \qquad \qquad
= z_2^{-1}    
\delta\left( \frac{z_1 - z_0}{z_2}\right)
Y \left( Y(u, z_0)v, z_2\right).  
\end{eqnarray} 
These axioms imply 
locality, skew-symmetry, associativity and commutativity conditions:
\begin{eqnarray}
&(z_{1}-z_{2})^N
Y(u,z_{1})Y(v,z_{2}) 
= 
(z_{1}-z_{2})^N
Y(v,z_{2})Y(u,z_{1}),&
\label{Local}
\\
\nonumber
\\
&Y(u,z)v = 
e^{zL(-1)}Y(v,-z)u,&
\label{skew}
\nonumber
\\
\nonumber
\\
\label{tornado}
&(z_{0}+z_{2})^N Y(u,z_{0}+z_{2})Y(v,z_{2})w = (z_{0}+z_{2})^N Y(Y(u,z_{0})v,z_{2})w,&
\label{Assoc}
\\
\nonumber
\\
&u(k)Y(v,z)- 
 Y(v,z)u(k)
= \sum\limits_{j\ge 0}  \left( k \atop j \right)
Y(u(j)v,z)z^{k-j},&
\label{Comm}
\end{eqnarray}
for $u$, $v$, $w\in V$ and integers  $N \gg 0 $.  
\end{definition}
For the conformal vector $\omega$ one has 
\begin{equation}
Y(\omega ,z)=\sum_{n\in {\mathbb Z}}L(n)z^{-n-2},  \label{Yomega}
\end{equation}
where $L(n)$ satisfies the Virasoro algebra for some central charge $C$ 
\begin{equation}
[\,  L(m),L(n)\, ]=(m-n)L(m+n)+\frac{C}{12}(m^{3}-m)\delta_{m,-n}{\rm Id}_V, 
\label{Virasoro}
\end{equation}
where ${\rm Id}_V$ is identity operator on $V$. 
Each vertex operator satisfies the translation property 
\begin{equation}
\partial_z Y(u,z)= Y\left(L(-1)u,z\right).  
\label{YL(-1)}
\end{equation}
The Virasoro operator $L(0)$ provides the ${\mathbb Z}$-grading with 
\[
L(0)u=ru, 
\]
 for 
$u\in V_{r}$, $r\in {\mathbb Z}$. 

For $v={\mathbf 1}_V$  
one has
\begin{equation}
\label{edinica}
Y({\mathbf 1}_V, z)={\rm Id}_V. 
\end{equation}
Note also that modes of homogeneous states are graded operators on $V$, i.e., for $v \in V_k$,  
\begin{eqnarray}
\label{modeaction1}
 v(n): V_m \rightarrow V_{m+k-n-1}.
\end{eqnarray}
 In particular, let us define the zero mode $o(v)$ of a state of weight $wt(v)=k$, i.e., $v \in V_k$,
as 
\begin{equation}
\label{zero_mode}
o(v) = v(wt (v) - 1), 
\end{equation}
 extending to $V$ additively. 
Similar definition is given for a vertex operator algebra module \cite{K}.  
\begin{definition}
Given a vertex operator algebra $V$, one defines the 
{adjoint vertex operator} with respect to $\alpha \in \C$,  by

\bigskip 
The adjoint operators defined by the mapping $z\mapsto 1/z$, i.e., 
with $\alpha=\epsilon$,  
\begin{eqnarray*}
&& Y_\epsilon^\dagger[v,z] = Y\left[\exp\left({\frac{z}{\epsilon}L[1]}\right)
\left(- \frac{\epsilon}{z^2}\right)^{L[0]}v, 
\frac{\epsilon}{z}\right].
\end{eqnarray*}
associated with the formal M\"obius map \cite{FHL}
\[
z \mapsto \frac{\alpha}{z}.  
\]
\end{definition}
\begin{definition}
An element $u\in V$ is called quasiprimary if 
\[
L(1)u=0. 
\]
\end{definition}
For quasiprimary $u$  
of weight $\wt(u)$ one has
\begin{align*}
u^\dagger(n) = (-1)^{\wt(u)} \alpha^{n+1-\wt(u)} u(2\wt(u)-n-2).
\end{align*}
\begin{definition}
A bilinear form 
\[
\langle . , . \rangle : V \times V\rightarrow \C, 
\]
 is called  {invariant} if \cite{FHL, Li}
\begin{align}
\langle Y(u,z)a,b \rangle = \langle a,Y^\dagger(u,z)b \rangle, 
\label{eq:bilform}
\end{align}
for all  
$a$, $b$, $u\in V$.
\end{definition}
Notice that the adjoint vertex operator $Y^\dagger(.,.)$ as well as the bilinear form $\langle. ,. \rangle$,  
 depend on $\alpha$.
 In terms of modes, we have 
\begin{align}
\langle u(n)a,b \rangle = \langle a,u^\dagger(n)b \rangle.
\label{eq:undag}
\end{align} 
Choosing $u=\omega$, and for $n=1$ implies 
\[
\langle L(0)a,b \rangle = \langle a,L(0)b \rangle.
\]
Thus, 
\[
\langle a, b \rangle=0, 
\]
 when $\wt(a) \neq \wt(b)$.  
\begin{definition}
A vertex operator algebra is called of strong-type 
if 
\[
V_{0} = \C\mathbf{1}_V, 
\]
and $V$ is simple and self-dual, i.e., $V$ is isomorphic to the dual module $V^\prime$ as a $V$-module.  
\end{definition}
It is proven in \cite{Li} 
that a strong-type vertex operator algebra $V$ has a unique invariant non-degenerate  
bilinear form up to normalization. 
This motivates 
\begin{definition}
\label{def:LiZ}
The form 
$\langle  ., .\rangle$
 on a strong-type vertex operator algebra $V$ is the unique invariant bilinear form $\langle . , . \rangle$  
normalized by 
\[
\langle \mathbf{1}_V, \mathbf{1}_V \rangle = 1.
\]
\end{definition}
Given a vertex operator algebra $(V,Y(.,.),\mathbf{1}_V,\omega)$, one can find an isomorphic vertex operator algebra 
 $(V,Y[.,.],\mathbf{1}_V,\omt)$ 
called \cite{Zhu} the 
{square-bracket} vertex operator algebra. 
 Both algebras have the same underlying vector space $V$, vacuum vector $\mathbf{1}_V$, and central charge. 
The vertex operator $Y[., .]$ is determined by 
\begin{align*}
Y[v,z] = \sum_{n\in \Z} v[n] z^{-n-1} = Y\left(q_z^{L(0)}v, q_z-1\right).
\end{align*} 
The new square-bracket conformal vector is 
\[
\omt = \omega - \frac{c}{24}\mathbf{1}_V,
\]
with the vertex operator
\[
Y[\omt,z] = \sum_{n\in \Z} L[n] z^{-n-2}. 
\]
The square-bracket Virasoro  operator mode 
$L[0]$ provides an alternative $\Z$--grading on $V$, i.e.,   
$\wt[v]=k$ if 
\[
L[0]v = kv, 
\]
 where $\wt[v]=\wt(v)$ for 
 primary $v$, and $L(n)v=0$ for all $ n>0$. 
 We can similarly define a square-bracket
bilinear form 
$\langle . , . \rangle_{\mathrm{sq}}$.
\section{Appendix: genus $g$ generalizations of elliptic functions}
In this Appendix we recall \cite{T2} genus $g$ generalizations of classical elliptic functions. 
\subsection{Classical elliptic functions}  
\label{subsect_Elliptic}
In this subsection we recall the classical elliptic functions and    
 modular forms \cite{Se, La}.
\begin{definition}
\label{def:Eisenstein}
For an integer $k\geq 2$, the {Eisenstein series}  is given by 
\begin{align*}
E_k(\tau) = E_k(q) = \delta_{n, even}  
\left(-\frac{B_{k}}{k!}+\frac{2}{(k-1)!}\sum_{n\geq 1}\sigma_{k-1}(n)q^{n}\right),  
\end{align*}
where $\tau\in \HH$,  $q=e^{2\pi i \tau}$,  
\[
\sigma_{k-1}(n) = \sum_{d\vert n} d^{k-1}, 
\]
 and $k-{\mathrm{th}}$ 
Bernoulli number  $B_k$.
\end{definition}
If $k\geq 4$ then $E_k(\tau)$ is a  modular form of weight $k$ on $\SL(2,\Z)$, 
while $E_2(\tau)$ is a quasi-modular form. 
\begin{definition}
\label{def:Pfunctions}
For integer $k\ge 1$, define elliptic functions $z\in \C$: 
\begin{eqnarray*}
&&
P_{k}(z,\tau)=\frac{(-1)^{k-1}}{(k-1)!}\partial_z^{k-1}P_{1}(z,\tau),
\nn
&&
 P_{1}(z,\tau)=\frac{1}{z}-\sum_{k\geq 2}E_{k}(\tau )z^{k-1}.
\end{eqnarray*}
\end{definition}
In particular 
\[
P_{2}(z, \tau) =\wp (z,\tau)+E_{2}(\tau ), 
\]
 for 
Weierstrass function $\wp (z,\tau)$ with periods $2\pi i$ and $2\pi i\tau$.
$P_{1}(z, \tau)$ is related to the quasi--periodic Weierstrass $\sigma$--function  with 
\[
P_{1}(z+2\pi i\tau, \tau)=P_{1}(z, \tau)-1. 
\]
\subsection{Genus two counterparts of Weierstrass functions}
\label{derivation}
In this section we recall the definition of genus two Weierstrass functions \cite{GT}.   
For $m$, $n\ge 1$, 
 we first define a number of infinite matrices and row and column vectors: 
\begin{align}
\Gamma(m,n) &= \delta_{m, -n+2p-2}, 
\nn
\Delta(m,n) &= \delta_{m, n+2p-2}. 
\label{eq:GamDelTh}
\end{align}
We also define the projection matrix
\begin{eqnarray}
\Pi =\Gamma^2=
\begin{bmatrix}
\mathbbm{1}_{2p-1} & 0 \\
0 & \ddots 
\end{bmatrix}
,
\label{eq:PiK}
\end{eqnarray}
where ${\rm Id}_{2p-3}$ denotes the $2p-3$ dimensional identity matrix and  
${\rm Id }_{-1}=0$.   
 Let $\Lambda_a $  for 
$a\in \{1,2\}$ 
be  the matrix with components
\begin{eqnarray}
\Lambda_a (m,n) &=& \Lambda_a (m,n;\tau_a,\epsilon) 
\nn
&=&
\epsilon^{(m+n)/2}(-1)^{n+1} \binom{m+n-1}{n}E_{m+n}(\tau_a).
\label{eq:Lambda}
\end{eqnarray}
Note that  
\begin{align}
\Lambda_a =SA_{a}S^{-1}, 
\label{eq:LamA}
\end{align}
for $A_a$ given by 
\begin{align*}
A_{a}= A_{a}(k,l,\tau _{a},\epsilon ) 
= \frac{ (-1)^{k+1}\epsilon ^{(k+l)/2}}{\sqrt{kl}} \frac{(k+l-1)!}{(k-1)!(l-1)!}E_{k+l}(\tau_a).
\end{align*}
introduce the infinite dimensional matrices 
 for $S$ a diagonal matrix with components
\begin{align}
S(m,n)=\sqrt{m}\delta_{mn}.
\label{eq:Sdef}
\end{align}
Let $\mathbb{R}(x) $ for $x\in \widehat{\Sigma}^{(1)}_a$ be the row vector with components
\begin{align}
\mathbb{R}(x;m) = \epsilon^{\frac{m}{2}} P_{m+1} (x, \tau_a). 
\label{eq:Rdef}
\end{align}
Let $\mathbb{X}_a$ for $a\in \{1,2\}$  be the column vector with components
\begin{eqnarray}
\label{eq:Xadef}
\mathbb{X}_1(m) &=& \mathbb{X}_1 \left(m; v_{n+1},  
 {\bf x}_n; B^{(2)} \right) 
\nn 
&=& \epsilon^{-m/2}\sum_{u\in V} 
\F^{(1)}_{V, k}\left( 
{\bf Y[x_k]}  
v[m]u; \tau_1  \right) \; 
\F^{(1)}_{V, n-k} \left( 
{\bf Y[x_{k+1, n}]  } \; 
\overline{u};\tau_2\right),  
\nn
\mathbb{X}_2(m) &=& 
\mathbb{X}_2\left(m; v_{n+1},
{\bf x}_n; B^{(2)} \right)  
\nn 
&=& \epsilon^{-m/2}\sum_{u\in V} 
\F^{(1)}_{V, k}\left( 
{\bf Y[x_k]} \;  
u;\tau_1\right)   
\F^{(1)}_{V, m-k}\left( 
{\bf Y[x_{n-k}]} 
v[m] \overline{u} ;\tau_2\right).
\end{eqnarray}
Introduce also 
$\mathbb{Q}(p; x)$  
 an infinite row vector defined  by
\begin{equation}
\mathbb{Q}(p; x) =\mathbb{R}(x) \Delta 
\left( \mathbbm{1} - \widetilde{\Lambda}_{\overline{a}}\widetilde{\Lambda}_a \right)^{-1},
\label{eq:Qdef}
\end{equation}
for $x\in\widehat{\Sigma}^{(1)}_a$.
Notice that  
\[
 \widetilde{\Lambda}_a=\Lambda_a\Delta.
\] 
On introduces 
\[
\mathbb{P}_{j+1}(x)=\frac{(-1)^j}{j!}\mathbb{P}_{1}(x),
\]
 for $x\in \widehat{\Sigma}^{(1)}_a$, and $j\ge 0$, is the column 
with components
\begin{align}
\mathbb{P}_{j+1}(x;m)=\epsilon^{\frac{m}{2}}\binom{m+j-1}{j}
\left( P_{j+m}(x,\tau_a) - \delta_{j 0}E_m(\tau_a)\right). 
\label{eq:P1jdef}
\end{align}
\begin{definition}
One defines  
\[
\mathcal{P}_{1}(p; x,y)=
\mathcal{P}_{1}(p; x,y; \tau_1, \tau_2, \epsilon),  
\]
   for $p \ge 1$ by
\begin{align*}
\mathcal{P}_{1}(p; x,y)  
=&
P_1(x-y,\tau_a)- P_1(x,\tau_a)  
\nn
- &
\mathbb{Q}(p; x) \widetilde{\Lambda}_{\abar} \,\mathbb{P}_{1} (y)
- 
(1-\delta_{p1})
\left(
\mathbb{Q}(p; x)\Lambda_{\abar}\right) (2p-2),
\end{align*}
for   $x$, $y\in\widehat{\Sigma}^{(1)}_a$, and  
\begin{eqnarray*}
\mathcal{P}_{1}(p; x,y)  
=&
(-1)^{p+1} \Big[
\mathbb{Q}(p; x) \mathbb{P}_{1}  (y) 
 + 
 (1-\delta_{p1}) \epsilon^{p-1} P_{2p-1}(x)
\nn
+&   
(1-\delta_{p1}) \left( 
\mathbb{Q}(p; x) \widetilde{\Lambda}_{\abar} \Lambda_a \right)  (2p-2) \Big], 
\end{eqnarray*}
for $x\in\widehat{\Sigma}^{(1)}_a$, 
$y\in\widehat{\Sigma}^{(1)}_{\abar}$. 
For $j> 0$, define 
\[
\mathcal{P}_{ j+1}(p; x,y)  
= \frac{1}{j!} \partial_y^j \left(
\mathcal{P}_{1}(p; x,y) \right),
\]
 i.e., 
\begin{align}
\mathcal{P}_{j+1}(p; x, y)=
\delta_{a, \bar{a} }P_{j+1}(x-y)+ (-1)^{j+1}.  
\mathbb{Q}(p; x) \left(\widetilde{\Lambda}_{\abar} \right)^{\delta_{a,\bar{a} } }\; \mathbb{P}_{j+1}(y),\; & 
\label{eq:PN21j}
\end{align}
for $x\in\widehat{\Sigma}^{(1)}_a$, $y\in\widehat{\Sigma}^{(1)}_{\bar{a}}$. 
One calls  
$%
\mathcal{P}_{j+1}(p; x,y)$ 
the genus two generalized Weierstrass functions.
\end{definition}
\subsection{Genus $g$ generalizations of elliptic functions}
\label{derivation}
For purposes of the formula \eqref{eq:ZhuGenusg} 
 we recall here certain definitions \cite{TW}. 
Define a column vector 
\[
X=(X_{a}(m)), 
\]
 indexed by $ m\ge 0$ and $ a\in\I$  with components 
\begin{align}\label{XamDef}
X_{a}(m)=\rho_{a}^{-\frac{m}{2}}\sum_{\bm{b}_{+}}Z^{(0)}(\ldots;u(m)b_{a}, w_{a};\ldots),
\end{align}
and a row vector 
\[
p(x)=(p_{a}(x,m)), 
\]
  for $m\ge 0, a\in\I$  with components 
\begin{align}\label{eq:pdef}
p_{a}(x,m)=\rho_{a}^{\frac{m}{2}}\del^{(0,m)}\psi_{p}^{(0)}(x,w_{a}).
\end{align}
Introduce the  
column vector 
\[
G=(G_{a}(m)), 
\]
 for $m\ge 0, a\in\I$, given by 
\begin{align*}
G=\sum_{k=1}^{n}\sum_{j\ge 0}\del_{k}^{(j)} \; q(y_{k})\; 
Z_{V}^{(g)}((u(j))_k {\bf x}_n),  
\end{align*}
where $q(y)=(q_{a}(y;m))$,  for $m\ge 0$, $a\in\I$, is a column vector with components 
\begin{align}
\label{eq:qdef}
q_{a}(y;m)=(-1)^{p}\rho_{a}^{\frac{m+1}{2}}\del^{(m,0)}\psi_{p}^{(0)}(w_{-a},y),
\end{align}
and 
\[
R=(R_{ab}(m,n)), 
\]
 for $m$, $n\ge 0$ and $a$, $b\in\I$ is a doubly indexed matrix with components 
\begin{align}
R_{ab}(m,n)=\begin{cases}(-1)^{p}\rho_{a}^{\frac{m+1}{2}}\rho_{b}^{\frac{n}{2}}\del^{(m,n)}\psi_{p}^{(0)}(w_{-a},w_{b}),&a\neq-b,\\ 
(-1)^{p}\rho_{a}^{\frac{m+n+1}{2}}\E_{m}^{n}(w_{-a}),&a=-b, 
\end{cases}
\label{eq:Rdef}
\end{align}
where
\begin{align}
\label{eq:Ejt}
\E_{m}^{n}(y)=\sum_{\ell=0}^{2p-2}\del^{(m)}f_{\ell}(y)\;\del^{(n)}y^{\ell}, 
\end{align}
\begin{align}\label{PsiDef}
\psi_{p}^{(0)}(x,y)=\frac{1}{x-y}+\sum_{\ell=0}^{2p-2}f_{\ell}(x)y^{\ell},
\end{align}
for {any} Laurent series $f_{\ell}(x)$ for $\ell=0,\ldots ,2p-2$.
Define the doubly indexed matrix $\Delta=(\Delta_{ab}(m,n))$ by
\begin{align}
\Delta_{ab}(m,n)=\delta_{m,n+2p-1}\delta_{ab}. 
\label {eq:Deltadef}
\end{align}
Denote by 
\[
\widetilde{R}=R\Delta,
\]
and the formal inverse $(I-\widetilde{R})^{-1}$ is given by 
\begin{align}
\label{eq:ImRinverse}
\left(I-\widetilde{R}\right)^{-1}=\sum_{k\ge 0}\widetilde{R}^{\,k}.
\end{align}
Define 
$\chi(x)=(\chi_{a}(x;\ell))$ and 
\[
o(u;\bm{v,y})=(o_{a}(u;\bm{v,y};\ell)), 
\]
 are 
{finite} row and column vectors indexed by 
$a\in\I$, $0\le \ell\le 2p-2$ with
\begin{align}
\label{eq:chiadef}
\chi_{a}(x;\ell)&=\rho_{a}^{-\frac{\ell}{2}}(p(x)+\widetilde{p}(x)(I-\widetilde{R})^{-1}R)_{a}(\ell),
\\
\label{LittleODef}
o_{a}(\ell)&=o_{a}(u;\bm{v,y};\ell)=\rho_{a}^{\frac{\ell}{2}}X_{a}(\ell),
\end{align}
and where 
\[
\widetilde{p}(x)=p(x)\Delta. 
\]
  $\psi_{p}(x,y)$ is defined by
\begin{align}
\label{eq:psilittleN}
\psi_{p}(x,y)=\psi_{p}^{(0)}(x,y)+\widetilde{p}(x)(I-\widetilde{R})^{-1}q(y).
\end{align}
For each $a \in\Ip$ we define a vector
\[
\theta_{a}(x)=(\theta_{a}(x;\ell) ), 
\]
 indexed by 
$0\le \ell\le  2p-2$ with components
\begin{align}\label{eq:thetadef}
\theta_{a}(x;\ell) = \chi_{a}(x;\ell)+(-1)^{p }\rho_{a}^{p-1-\ell}\chi_{-a}(x;2p-2-\ell).
\end{align}
Now define the following vectors of formal differential forms
\begin{align}
\label{eq:ThetaPQdef}
 P(x) =p(x) \; dx^{p},
\nn
 Q(y)=q(y)\; dy^{1-p},
\end{align}
with 
\[
\widetilde{P}(x)=P(x)\Delta.
\]
 Then with 
\begin{equation}
\label{psih}
\Psi_{p} (x,y) =\psi_{p}(x,y) \;dx^{p}\; dy^{1-p}, 
\end{equation}
 we have
\begin{align}\label{GenusgPsiDef}
\Psi_{p}(x,y)=\Psi_{p}^{(0)}(x,y)+\widetilde{P}(x)(I-\widetilde{R})^{-1}Q(y).
\end{align}
Defining   
\begin{equation}
\label{thetanew}
\Theta_{a}(x;\ell) =\theta_{a}(x;\ell)\; dx^{p}, 
\end{equation}
  and 
\begin{equation}
\label{oat}
O_{a}(u; \bm{v,y};\ell) = o_{a}(u; \bm{v,y};\ell) \; \bm{dy^{\wt(v)}}, 
\end{equation}
\begin{remark}
	\label{rem:MainTheorem}
The $\Theta_{a}(x)$, and $\Psi_{p}(x,y )$ coefficients depend 
on $p=\wt(u)$ but  are otherwise independent of the vertex operator algebra $V$. 
\end{remark}

\end{document}